\documentclass{article}
\usepackage{amsfonts}
\usepackage{amsmath,amsthm}

\usepackage[all]{xy}
\usepackage{graphicx}

\usepackage{amssymb}
\usepackage{latexsym}

\DeclareMathAlphabet\EuFrak{U}{euf}{m}{n}	
\SetMathAlphabet\EuFrak{bold}{U}{euf}{b}{n}	

\newcommand{\ra}{\rightarrow}

\newcommand{\hra}{\hookrightarrow}

\newcommand{\ovl}{\overline}

\newcommand{\wa}{\widehat}

\newcommand{\sC}{{\it C*}-}
\newcommand{\bC}{{\mathbb C}}
\newcommand{\bR}{{\mathbb R}}

\newcommand{\bT}{{\mathbb T}}

\newcommand{\bZ}{{\mathbb Z}}
\newcommand{\bM}{{\mathbb M}}
\newcommand{\bN}{{\mathbb N}}

\newcommand{\ud}{{{\mathbb U}(d)}}

\newcommand{\eps}{\varepsilon}

\newcommand{\mA}{\mathcal A}
\newcommand{\mB}{\mathcal B}
\newcommand{\mC}{\mathcal C}

\newcommand{\mE}{\mathcal E}
\newcommand{\mF}{\mathcal F}

\newcommand{\mH}{\mathcal H}

\newcommand{\mK}{\mathcal K}
\newcommand{\mL}{\mathcal L}
\newcommand{\mM}{\mathcal M}
\newcommand{\mN}{\mathcal N}
\newcommand{\mO}{\mathcal O}

\newcommand{\mR}{\mathcal R}
\newcommand{\mT}{\mathcal T}

\newcommand{\mU}{\mathcal U}
\newcommand{\mV}{\mathcal V}
\newcommand{\mZ}{\mathcal Z}

\newcommand{\ii}{\iota,\iota}

\newcommand{\mrs}{\mM^r,\mM^s}
\newcommand{\hrs}{\mH^r,\mH^s}
\newcommand{\ers}{\mE^r,\mE^s}

\newcommand{\wE}{\wa {\mathcal E}}
\newcommand{\wH}{\wa {\mathcal H}}

\newcommand{\zro}{\mZ^{\rho}}

\newcommand{\coe}{\mO_{\mE}}
\newcommand{\coh}{\mO_{\mH}}
\newcommand{\com}{\mO_{\mM}}

\newcommand{\toe}{\mT_{\mE}}
\newcommand{\toh}{\mT_{\mH}}

\newcommand{\doe}{{\widetilde {\mO}}_{\mE}}
\newcommand{\doh}{{\widetilde {\mO}}_{\mH}}
\newcommand{\dom}{{\widetilde {\mO}}_{\mM}}

\newtheorem{thm}{Theorem}[section]
\newtheorem{cor}[thm]{Corollary}
\newtheorem{lem}[thm]{Lemma}
\newtheorem{prop}[thm]{Proposition}

\newtheorem{defn}[thm]{Definition}

\theoremstyle{definition}
\newtheorem{ex}{Example}[section]

\theoremstyle{remark}
\newtheorem{rem}{Remark}[section]

\numberwithin{equation}{section}

\begin{document}

\author{{\sf Ezio Vasselli}
                         \\{\it Dipartimento di Matematica}
                         \\{\it Universit\`a di Roma ``Tor~Vergata''}
			 \\{\it Via della Ricerca Scientifica, snc - 00133 Roma - Italy }
                         \\{\sf vasselli@mat.uniroma2.it}}

\title{ The \sC algebra of a vector bundle \\and\\
        fields of Cuntz algebras}
\maketitle

\begin{abstract}
We study the Pimsner algebra associated with the module of continuous sections of a Hilbert bundle, and prove that it is a continuous bundle of Cuntz algebras. We discuss the role of such Pimsner algebras w.r.t. the notion of inner endomorphism. Furthermore, we study bundles of Cuntz algebras carrying a global circle action, and assign to them a class in the representable $KK$-group of the zero-grade bundle. We explicitly compute such class for the Pimsner algebra of a vector bundle.

\bigskip

\noindent {\bf AMS Subj. Class.:} 46L05, 46L80.

\noindent {\bf Keywords:} Pimsner algebra; Cuntz algebra; Continuous Field; Vector Bundle.

\end{abstract}

\section{Introduction.}
\label{intro}

The Pimsner algebras (also called Cuntz-Pimsner algebras, or Cuntz-Krieger-Pimsner algebras) were introduced in \cite{Pim93}, as a natural generalization of the Cuntz-Krieger algebras (\cite{Cun77,CK80}) and crossed products by a single automorphism. They also include crossed products by partial automorphisms and by Hilbert bimodules in the sense of \cite{Exe94,AEE95}. A universal construction, generalizing the one of the Pimsner algebra, has been made by Doplicher and Roberts in the framework of tensor categories (\cite[\S 3]{DR89}). The so-constructed \sC algebras are called DR-algebras (see \cite{DR89,LR97,Vas}). Studies of the Pimsner algebra from this point of view are made in \cite{DPZ97,KPW98,KPW01}.

Let $X$ be a locally compact Hausdorff space. In the present paper we proceed to the study of the Pimsner (and DR) algebra associated with the module of continuous sections of a Hilbert bundle $\mH \ra X$, started in \cite{Vas} in a categorical framework (in the sequel, we will more concisely say {\em the Pimsner algebra of the Hilbert bundle}). The interest in such algebras arises for their role in the context of crossed products by endomorphisms, and duality for non-compact groups (see \cite{Vas03}). Anyway, they are an interesting class also from the viewpoint of the classification of \sC algebras by $KK$-theoretical invariants, in the spirit of \cite{Kir00,BK02}. We will develop here the basic properties of such Pimsner algebras, and give some applications for Cuntz algebra bundles; more detailed $K$-theoretical and classification questions will be approached in future papers. The present work is organized as follows:

\bigskip

In Section \ref{preliminaries} we give some preliminary results about continuous bundles of \sC algebras, Pimsner and DR-algebras, representable $KK$-theory.

\bigskip

In Section \ref{bimod_alg} we discuss some basic properties of Hilbert bimodules in \sC algebras, and their relationship with endomorphisms. In this setting, we give a generalized notion of inner endomorphism (Def.\ref{inner}), and give some motivating examples.

\bigskip

In Section \ref{basics} we investigate the basic properties of the Pimsner algebra $\coh$ of a Hilbert bundle $\mH \ra X$. We prove that $\coh$ is a continuous bundle, having as fibre the Cuntz algebra $\mO_d$ (Prop.\ref{coe_cf}). In the case in which $\mE \ra X$ is a vector bundle, we prove that the DR-algebra $\doe$ is the multiplier algebra of $\coe$ (Prop.\ref{coe_doe}), and relate $\doe$ with the Pimsner algebra of a 'pullback bundle' $\beta \mE \ra \beta X$ over the Stone-Cech compactification (Prop.\ref{sigma_inner}). We also consider the canonical endomorphism over $\coe$ (resp. $\doe$), and discuss the connection with inner endomorphisms (Prop.\ref{sigma_uni}). When $X$ is compact, a simple description of $\coe$ can be made in terms of generators and relations (Thm.\ref{thm_genrel}).

\bigskip

In Section \ref{homotopy} we study graded $\mO_d$-bundles (i.e., locally trivial continuous bundles of Cuntz algebras carrying a global circle action). We give an isomorphism condition in terms of a suitable Hilbert bimodule (the grade one component, see Prop.\ref{thm_f1}), which allows to assign a $KK$-theoretical invariant to the given graded $\mO_d$-bundle, belonging to the representable Kasparov group of the zero grade \sC algebra (Eq.(\ref{def_delta_1})). Such invariant is used in a Pimsner-Voiculescu exact sequence for the $KK$-theory of the $\mO_d$-bundle (Cor.\ref{cor_ex_seq}), and is computed for the Pimsner algebra of a vector bundle (Thm.\ref{delta1_coe}). As an application, we prove that a stable graded isomorphism at level of the Pimsner algebra implies an equivalence in terms of (representable) $K$-theory of the underlying vector bundles (Prop.\ref{prop_k}). Viceversa, if the base space is a finite $CW$-complex, the equivalence in $K$-theory of the vector bundles implies a stable isomorphism of the Pimsner algebras (Prop.\ref{prop_cw}).

\bigskip

In the final Section \ref{final} we briefly discuss further developement about classification questions, relationships with the duality theory for non-compact groups and projective multiresolution analysis in the sense of \cite{PR03}.

\section{Preliminaries.}
\label{preliminaries}

\subsection{Fields, Bundles.}
\label{fields}

The notion of continuous field of \sC algebras (resp. Banach spaces, Hilbert spaces) is well known in literature; for basic notions and terminology we will refere to  \cite[\S 10]{Dix}. 

\bigskip

In particular, we consider continuous fields of Hilbert (Banach) spaces having finite dimensional fibres; we call {\em Hilbert (Banach) bundles} the corresponding {\em espaces fibr\'es vectoriels} in the sense of \cite[I.2]{DD63}, according to the terminology of Dupr\'e \cite{Dup74,Dup76}. We recall that the {\em order} of a Hilbert (Banach) bundle is the cardinality of the set of isomorphism classes of the fibres; if a Hilbert bundle has order one, then it is a vector bundle (i.e. the Hilbert bundle is locally trivial, see \cite[Prop.2.3]{Dup74}).

Let $X$ be a locally compact Hausdorff space. We will adopt the notation $p : \mH \ra X$ for a Hilbert (Banach) bundle, while $\wH$ will denote the Hilbert (Banach) $C_0(X)$-bimodule of continuous, vanishing at infinity sections of $\mH$. In the present work we will always make the standard assumption that $\mH$ is full, i.e. for every $v \in \mH$, $p(v) = x$, there is $\psi \in \wH$ such that $v = \psi (x)$. We will call the {\em rank} of $\mH$ the map $d : X \ra \bN$ assigning to $x \in X$ the dimension of the fibre of $\mH$ over $x$; if $\mH$ is a vector bundle, then $d$ is continuous.

In particular, a vector bundle will be denoted by the usual notation $p : \mE \ra X$. For basic properties and terminology about vector bundles we refere to \cite{Ati,Kar}.

It is well-known from the Serre-Swan theorem that in the case in which the base space $X$ is compact, the category of vector bundles is equivalent to the one of projective, full, finitely generated Hilbert $C(X)$-modules. If we pass to the category of Hilbert bundles, then some properties are lost at level of the corresponding module $\wH$ of continuous sections; for example, $\wH$ could be not full, but this case can be easily avoided by regarding $\wH$ as a module over the \sC algebra of scalar products (which is an ideal in $C(X)$), and by rescaling the base space to be the corresponding open subset of $X$. {\em In the present paper, we will consider Hilbert bundles having full module of continuous sections}. Note that this property should not be confused with the fact that the bundle is full.

For the notions of morphism, tensor product and dual (conjugate) of a Hilbert bundle, we will refere to \cite[II.18]{DD63}, \cite[\S 4]{Dup74}.

We will also consider the following notion (\cite{KW95}): a {\em continuous bundle} of \sC algebras over $X$ is a \sC algebra $\mA$, equipped with a faithful family of epimorphisms $\left\{ \pi_x : \mA \ra \mA_x \right\}_{x \in X}$ such that, for every $a \in \mA$, the {\em norm function} $\left\{ x \mapsto \left\| \pi_x (a) \right\|  \right\}$ belongs to $C_0(X)$; furthermore, $\mA$ is required to be a nondegenerate $C_0(X)$-module w.r.t. pointwise multiplication $f , a \mapsto \left\{ f(x) \cdot \pi_x (a) \right\}$, $f \in C_0(X)$. The \sC algebras $\mA_x$ are called the {\em fibres} of $\mA$ over $x$. Note that $\mA$ may be regarded as the \sC algebra of a continuous field in the sense of \cite[\S 10]{Dix}; thus, $\mA$ is the \sC algebra of continuous, vanishing at infinity sections of a topological bundle $\wa \mA \ra X$, defined as the set $\wa \mA := \bigcup_x \mA_x$ endowed with the topology \cite[I.2]{DD63}.

We can relax the properties required in the definition of continuous bundle, by assuming that the norm function is upper semicontinuous. It is proven in \cite{Bla96,Nil96} that the category of such {\em upper semicontinuous bundles} coincides with the category of $C_0(X)$-{\em algebras}. A $C_0(X)$-algebra is a \sC algebra $\mA$, equipped with a nondegenerate morphism from $C_0(X)$ into the centre of the multiplier algebra $M(\mA)$. In the sequel we will identify $C_0(X)$ with the image in $M(\mA)$. 

\bigskip

We will denote by $\otimes_X$ the minimal $C_0(X)$-algebra tensor product (\cite{Bla95}).

\bigskip

Let $\mA$ be a \sC algebra. In order for a more concise terminology, we will call $\mA${\em - bundle} a locally trivial continuous bundle with fibre $\mA$.

\bigskip

Let $X$ be a locally compact Hausdorff space, $C_b(X)$ the \sC algebra of continuous bounded functions over $X$. We denote by $\beta X$ the Stone-Cech compactification of $X$, so that $C(\beta X) \simeq C_b (X)$. Let $\mA$ be a unital \sC algebra, $\mF$ an $\mA$-bundle over $X$. We denote by $C_b (X , \wa \mF)$ the \sC algebra whose elements are the bounded continuous sections of the topological bundle $\wa \mF \ra X$ associated with $\mF$. If $M(\mF)$ is the multiplier algebra, then by \cite[Thm.3.3]{APT73} there is an isomorphism $M(\mF)  \simeq  C_b (X , \wa \mF)$; note that in order to get such an isomorphism, both the properties $\mA \ni 1$, $\mF$ locally trivial, are required.

\bigskip

In order for more concise notations, for every \sC algebra $\mA$ we will denote 
\begin{equation}
\label{not_tp}
X \mA := C_0 (X , \mA) \ \ , \ \ I \mA := C([0,1] , \mA) \ \ .
\end{equation}

\subsection{Pimsner algebras vs. DR-algebras.}
\label{pimsner}

We give in the present section a brief account about the Pimsner algebra associated with a Hilbert bimodule (\cite{Pim93}), as it was constructed in \cite[\S 3]{DPZ97} in the abstract setting of (semi)tensor \sC categories.

Let $\mA$ be a \sC algebra, $\mM$ a Hilbert $\mA$-bimodule, which is assumed to be full and with non-degenerate left $\mA$-module action; the $\mA$-valued scalar product is denoted by $\left \langle \cdot , \cdot \right \rangle$. We consider the Banach $\mA$-bimodules $(\mrs)$ of bounded right $\mA$-module operators from the tensor power $\mM^r$ into $\mM^s$, with $r,s \in \bN$ (we are considering the {\em internal} tensor product of Hilbert bimodules). For $r = 0$, we define $\mM^0 := \mA$ and $(\mA , \mA) := M(\mA)$, where $M(\mA)$ is the multiplier algebra. Let $1 \in (\mM , \mM)$ be the identity map on $\mM$; we introduce an injective map $i : (\mrs) \hra (\mM^{r+1}, \mM^{s+1})$, by defining $i(t) := t \otimes 1$. Then, we consider the inductive limit 
\[  \dom^k := \cdots \hra (\mrs) \hra (\mM^{r+1}, \mM^{s+1}) \hra \cdots  \ \ ,  \]
\noindent where $k := s - r \in \bZ$. Let $t \in \dom^h$, $t' \in \dom^k$; we can always assume to find representatives $t \in (\mM^{r}, \mM^{r+h})$, $t' \in (\mM^{r+h}, \mM^{r+h+k})$ in the inductive limits $\dom^h$, $\dom^k$. We define a product, by applying the composition $t' \cdot t \in (\mM^r , \mM^{r+h+k}) \subset \dom^{h+k}$. In the same way, the operation of assigning the adjoint operator induces an involution $* : \dom^k \ra \dom^{-k}$. We define 
\[  \ ^0\dom := {\sum_k}^\oplus \dom^k \ \ . \]
\noindent The above considerations imply that $^0\dom$ is a $*$-algebra. It can be proved that there exists a unique \sC norm over $^0\dom$ such that the {\em circle action} 
\begin{equation}
\label{def_circ_act}
\wa z  \ \left( {\sum}_k t_k \right)  \ := \ {\sum}_k z^k t_k  \ \ ,
\end{equation}
\noindent $z \in \bT$, $t_k \in \dom^k$, extends to an isometric action by automorphisms. We define the DR-algebra $\dom$ as the closure w.r.t. such \sC norm.

The Pimsner algebra $\com$ is defined as the \sC subalgebra of $\dom$ generated by the spaces of compact right $\mA$-module operators $\mK (\mrs) \subseteq (\mrs)$. We define $\mK (\mA , \mA) := \mA$, and note that there is a natural identification $\mM = \mK (\mA , \mM)$. In the same way as above, we define
\[  
\com^k := \bigcup_{r - s = k} \mK (\mM^r , \mM^s) \subseteq \dom^k  \ \ , \ \
^0\com := {\sum_k}^\oplus \com^k  \subseteq \ ^0\dom \ \ ,
\]
\noindent so that the circle action is well-defined also over $\com$. Of course, there is a monomorphism from $\com$ into $\dom$. If $\mM$ is projective, full and finitely generated, the DR-algebra $\dom$ coincides with $\com$.

\bigskip

\begin{rem}
\label{rem_theta}
Let $\psi , \psi' \in \mM^r$, $\varphi \in \mM^s$, and $\theta_{\varphi , \psi} \in \mK (\mrs) : \theta_{\varphi , \psi} (\psi') := \varphi \left \langle \psi , \psi' \right \rangle$. Then the following relations hold in the Pimsner (DR) algebra (see\cite[\S 3]{Pim93}): $\varphi \psi^* =  \theta_{\varphi , \psi} \ $, $ \ \psi^* \psi' = \left \langle \psi , \psi' \right \rangle$.
\end{rem}

\bigskip

Let $\mA = \bC$, and $\mM$ be an infinite dimensional Hilbert space; then the Pimsner algebra $\com$ is strictly contained in $\dom$ (see \cite{CDPR94}). If $\mM = \bC^d$, $d \in \bN$, $d > 1$, then $\com = \dom$ is the Cuntz algebra $\mO_d$ (\cite{Cun77}). The spectral subspaces w.r.t. the circle action over $\mO_d$ will be denoted by $\mO_d^k$, $k \in \bZ$; thus, $\mO_d^0$ is the UHF algebra $\otimes^\infty \bM_d$. In particular, if $\mM = \bC$ then $\dom$ is the \sC algebra $C(S^1)$ of continuous functions over the circle.

\bigskip

The Pimsner (DR) algebra has useful functoriality properties. In order to expose them, we consider some constructions in the category of Hilbert bimodules:

\begin{defn}
\label{2morphism}
Let $\mA$, $\mA'$ be a \sC algebras, $\mM$ a Hilbert $\mA$-bimodule, $\mN$ a Hilbert $\mA'$-bimodule. A {\bf covariant morphism} from $\mM$ into $\mN$ is a pair $(\beta , \eta)$, where $\beta : \mM \ra \mN$ is a Banach space map, $\eta : \mA \ra \mA'$ is a \sC algebra morphism, and the following properties are satisfied for $a \in \mA$, $\psi , \psi' \in \mM$:
\[  
\beta (a \psi) = \eta (a) \beta (\psi) \ \ , \ \ 
\beta (\psi a) = \beta (\psi) \eta (a) \ \ , \ \
\left \langle \beta (\psi) , \beta(\psi') \right \rangle  = 
      \eta \left \langle \psi , \psi' \right \rangle \  .
\]
\end{defn}

\begin{ex}
Let $G$ be a topological group, $\mA$ a \sC algebra, $\mM$ a $G$-Hilbert $\mA$-bimodule in the sense of \cite{Kas80}. Then $G$ acts on $\mM$ by covariant automorphisms (see also \cite[Rem.4.10(2)]{Pim93}).
\end{ex}

\begin{ex}
\label{ex_rest}
Let $X$ be a locally compact Hausdorff space, $Y \subset X$ a closed set, $\mM$ a countably generated Hilbert $C_0(X)$-bimodule such that $f \psi = \psi f$, $\psi \in \mM$, $f \in C_0(X)$. Then $\mM$ is a direct summand of $l^2 (C_0(X))$ by the Kasparov stabilization theorem. Let us denote by $\pi : C_0(X) \ra C_b (Y)$ the natural projection. Then a covariant morphism $\pi_*$ from $l^2 (C_0(X))$ onto $l^2(C_b(Y))$ is naturally induced, by defining $(\pi_* (F))_n := (\pi(F_n))$, where $F := (F_n) \in l^2(C_0(X))$. $\pi_*$ clearly restricts on a covariant morphism from $\mM$ into $l^2(C_b(Y))$.
\end{ex}

\bigskip

Let $\mA$ be a unital \sC algebra. Finitely generated, projective Hilbert $\mA$-bimodules can be described by {\em amplimorphisms}, i.e. \sC algebra morphisms $\phi : \mA \ra \mA \otimes \bM_n$. The Hilbert $\mA$-bimodule associated with $\phi$ is given by $\mM_\phi := \phi (1) (\oplus^n \mA)$, endowed with the obvious right $\mA$-module action and scalar product, and with left $\mA$-module action given by $a \cdot \psi := \phi (a) \psi$, $\psi \in \mM_\phi \subseteq \oplus^n \mA$. Note that if $\phi$ is unital, then $\mM_\phi$ is free as a right Hilbert $\mA$-module. Viceversa, given a projective Hilbert $\mA$-bimodule $\mM$ generated by a finite set $\left\{ \psi_l \right\}_{l=1}^n$, then an amplimorphism associated with $\mM$ is given by $\phi_\mM (a) := \left \langle  \psi_l , a \psi_m \right \rangle_{l,m} \in \mA \otimes \bM_n$. We leave to the reader the easy verification that $\mM$ is isomorphic to $\mM_{\phi_\mM}$ as a Hilbert $\mA$-bimodule.

\bigskip

\begin{ex}
\label{ex_pullback}
Let $\alpha : \mA \ra \mB$ be a \sC algebra isomorphism, $\mM$ a Hilbert $\mA$-bimodule. We introduce a Hilbert $\mB$-bimodule $\mM_\alpha$, defined as the set $\mM_\alpha := \left\{ \underline \psi , \psi \in \mM \right\}$ endowed with the vector space structure induced by $\mM$. The Hilbert $\mB$-bimodule structure is defined as follows:
\[
b \underline \psi := \underline {\alpha^{-1} (b) \psi} \ \ , \ \
\underline \psi b := \underline {\psi \alpha^{-1} (b)} \ \ , \ \
\left \langle \underline \psi , \underline \psi' \right \rangle := 
       \alpha \left \langle \psi , \psi' \right \rangle
\]
\noindent Let now $\beta (\psi) := \underline \psi$, $\psi \in \mM$. It is clear that the pair $(\beta , \alpha)$ is a covariant isomorphism from $\mM$ onto $\mM_\alpha$. We call $\mM_\alpha$ the {\em pullback bimodule} of $\mM$. Note that if $\mA$ is unital, $\mM$ is projective and $\phi : \mA \ra \mA \otimes \bM_n$ is an amplimorphism associated with $\mM$, then $\phi_\alpha := (\alpha \otimes id_n) \circ \phi \circ \alpha^{-1}$ is an amplimorphism associated with $\mM_\alpha$. 
\end{ex}

From the construction of the previous example, we deduce

\begin{lem}
\label{lem_pullback}
Let $\alpha : \mA \ra \mB$ be a \sC algebra isomorphism, $\mM$ a Hilbert $\mA$-bimodule, $\mN$ a Hilbert $\mB$-bimodule, $(\beta , \alpha) : \mM \ra \mN$ a covariant isomorphism. Then, the pullback bimodule $\mM_\alpha$ is isomorphic to $\mN$ as a Hilbert $\mB$-bimodule.
\end{lem}

Let $\mA$, $\mB$ be \sC algebras endowed with $\bT$-actions by automorphisms. A morphism $\phi : \mA \ra \mB$ intertwining the $\bT$-actions is said {\em graded}: in such a case, we will use the notation
\[  
\phi : ( \mA , \bT) \ra ( \mB , \bT ) \ \ .
\]

The following results are immediate consequence of the universality of the Pimsner (DR) algebra, so that the proofs are sketched or omitted (see \cite[\S 3]{Pim93}, \cite[\S 3]{DPZ97} for details).

\begin{lem}
\label{lem_funct}
Let $(\beta , \eta) : \mM \ra \mN$ be a covariant morphism. Then, there is a unique \sC algebra morphism $\beta_* : ( \com , \bT ) \ra ( \mO_{\mN} , \bT )$ extending $\beta$.
\end{lem}

\begin{proof}
For every $r \in \bN$, there is a covariant morphism $(\beta_r , \eta) : \mM^r \ra \mN^r$, $\beta_r (\psi_1 \otimes_\mA \cdots \otimes_\mA \psi_r) := \beta (\psi_1) \otimes_{\mA'} \cdots \otimes_{\mA'} \beta(\psi_r)$.  Let $\theta_{\psi' , \psi} \in \mK (\mrs)$, with $\psi \in \mM^r$, $\psi' \in \mM^s$. We define
\[  
\beta_* (\theta_{\psi' , \psi}) := \beta_s (\psi') \ \beta_r (\psi)^* \ \ , \ \ 
\beta_* (\theta_{\psi' , \psi}) \in \mK (\mN^r , \mN^s)  \ \ .
\]
\end{proof}

\begin{cor}
\label{u_act}
Let $\mM$ be a Hilbert $\mA$-bimodule, ${\bf U} \mM$ the group of unitaries in $\mB (\mM , \mM)$. Then, there is an action by $\mA$-bimodule automorphisms ${\bf U} \mM \ra {\bf aut} ( \dom , \bT )$, ${\bf U}\mM \ni u \mapsto \wa u :$ $\wa u (t) := u^{{\otimes}^s} \cdot t \cdot {u^*}^{{\otimes}^r}$, $t \in (\mrs)$.
\end{cor}

\begin{cor}
Let $G$ be a topological group, $\mM$ a $G$-Hilbert $\mA$-bimodule. Then $G$ acts continuously by automorphisms on $\dom$.
\end{cor}

\subsection{$C_0(X)$-Hilbert bimodules and $KK(X ; -,-)$-groups.}
\label{kk_groups}

For basic notions and terminology about $KK$-theory, we refere to \cite[\S 17]{Bla}.

\bigskip

Let $X$ be a locally compact Hausdorff space, $\mA$, $\mB$ $C_0(X)$-algebras. A $C_0(X)$-{\em Hilbert} $\mA$-$\mB$-{\em bimodule} is a Hilbert $\mA$-$\mB$-bimodule $\mM$ such that $(af) \psi b = a \psi (f b)$ for every $f \in C_0(X)$, $\psi \in \mM$, $a \in \mA$, $b \in \mB$.

\begin{ex}
Let $\mA$ be a \sC algebra, $X$ a locally compact Hausdorff space, $\mE \ra X$ a vector $\mA$-bundle in the sense of Mishchenko (\cite{MF80}). Then, the module of continuous, vanishing at infinity sections of $\mE$ has an obvious structure of $C_0(X)$-Hilbert $(X \mA)$-bimodule.
\end{ex}

\bigskip

Let now $X$ be a $\sigma$-compact metrisable space, and $\mA$, $\mB$ separable, $\bZ_2$-graded $C_0(X)$-algebras. A {\em Kasparov module} is a pair $(\mM , F)$, where $\mM$ is a countably generated, $\bZ_2$-graded $C_0(X)$-Hilbert $\mA$-$\mB$-bimodule, and $F = F^* \in (\mM , \mM)$ is an operator with degree $1$ such that $[F,a] , a(F^2 - 1) \in \mK (\mM , \mM)$ for every $a \in \mA$. The {\em representable} $KK${\em -theory group} $KK ( X ; \mA , \mB )$ is constructed in the same way as the usual $KK$-group, by endowing the set of homotopy classes of Kasparov modules with the operation of direct sum (see \cite[\S 2.19]{Kas88} for details: note that the notation $\mR KK (X ; \mA , \mB)$ is used instead of $KK (X ; \mA , \mB)$). 

In particular, $KK ( X ; \mA , \mA )$ has also been called the {\em central} $K$-{\em theory} of $\mA$ in \cite{PT00}, where an $E$-theory approach is defined in the $K(X)$-nuclear case. It is clear that when $X = \bullet$ reduces to a single point, then $KK ( \bullet , \mA , \mB )$ is the usual $KK$-theory group $KK_0 (\mA , \mB)$. A generalization of $KK (X ; \mA , \mB)$ to the case in which $X$ is a $T_0$ space has been made in \cite{Kir00}.

In the sequel, if not specified, every \sC algebra ($C_0(X)$-algebra) will be considered endowed with the trivial $\bZ_2$-grading, and the same assumption will be made for ($C_0(X)$)-Hilbert bimodules. We now introduce the notation
\[  RK^0 (X) := KK ( X ; C_0(X) , C_0(X) ) \ \ , \]
\noindent according to \cite[2.19]{Kas88}. If $X$ is compact, it is verified that $RK^0(X)$ coincides with the usual $K$-theory group $K^0(X)$.

\bigskip

The Kasparov product \cite[\S 2.21]{Kas88} induces a natural ring structure on $KK ( X ;\newline \mA , \mA )$, having as identity the class $[1]_\mA := [(\mA , 0)]$ associated with the obvious bimodule structure on $\mA$. If $n \in \bN$, we denote 
\begin{equation}
\label{eq_nkk}
[n] := [\oplus^n ( \mA , 0 ) ] \in KK( X ; \mA , \mA )  \ .
\end{equation}

\bigskip

Let $\mM$ be a $C_0(X)$-Hilbert bimodule, $\mC$ a $C_0(X)$-algebra. The algebraic tensor product $\mM \odot_{C_0(X)} \mC$ with coefficients in $C_0(X)$ is endowed with the following Hilbert $(\mA \otimes_X \mC)$-$(\mB \otimes_X \mC)$-bimodule structure
\begin{equation}
\label{def_ea}
\left\{
\begin{array}{ll}
(b \otimes c) , (\psi \otimes c') \mapsto (b \psi) \otimes (cc') \\

(\psi \otimes c') , (a \otimes c)  \mapsto (\psi a) \otimes (c'c) \\

\left \langle \psi \otimes c , \psi' \otimes c' \right \rangle := 
\left \langle \psi , \psi' \right \rangle \otimes_X (c^* c') \ \in \ \mB \otimes_X \mC

\end{array}
\right.
\end{equation}
\noindent where $\psi , \psi' \in \mM$, $a \in \mA$, $b \in \mB$, $c , c' \in \mC$. We denote by $\mM \otimes_X \mC$ the corresponding completition w.r.t. the $(\mB \otimes_X \mC)$-valued scalar product. it is clear that $\mM \otimes_X \mC$ is a $C_0(X)$-Hilbert $(\mA \otimes_X \mC)$-$(\mB \otimes_X \mC)$-bimodule.

\bigskip

As a consequence of the above construction, $KK ( X ; \mA , \mB )$ is endowed with a natural sheaf structure w.r.t. open subsets $Y \subseteq X$. In fact, if $Y$ is an open subset of $X$, then $C_0(Y)$ is a $C_0(X)$-algebra, endowed with the restriction morphism $\pi_Y : C_0(X) \ra M(C_0(Y))$, $\pi_Y (f) := f |_Y$. If $\mM$ is a $C_0(X)$-Hilbert $\mA$-$\mB$-bimodule, then $\mM \otimes_X C_0(Y)$ has a natural structure of $C_0(Y)$-Hilbert $\mA_Y$-bimodule, where $\mA_Y := \mA \otimes_X C_0(Y)$ (resp. $\mB_Y$) is the restriction of $\mA$ (resp. $\mB$) over $Y$ as an upper semicontinuous bundle. Thus, we obtain a map $KK ( X ; \mA , \mB ) \ra KK ( Y ; \mA_Y , \mB_Y )$.

\bigskip

\begin{rem}
\label{rem_cla}
Let $\mM$ be a countably generated $C_0(X)$-Hilbert $\mA$-$\mB$-bimodule such that $\mA$ acts on the left by elements of $\mK (\mM , \mM)$. Then, the pair $(\mM , 0)$ is a Kasparov module, and the corresponding cycle $[\mM] \in KK ( X ; \mA , \mB )$ is defined.
\end{rem}

\bigskip

Let $\mE \ra X$ be a vector bundle (with finite rank); we denote by $\wE$ the Hilbert $C_0(X)$-bimodule of continuous, vanishing at infinity section of $\mE$, endowed with the left $C_0(X)$-action coinciding with the right one.

\begin{lem}
\label{lem_vb_kk}
Let $X$ be a $\sigma$-compact Hausdorff space, $\mE \ra X$ a vector bundle. Then, $C_0(X)$ acts on the left over $\wE$ by elements of $\mK ( \wE , \wE )$, and the pair $(\wE , 0)$ is a Kasparov module with class $[\mE] := [(\wE , 0)] \in RK^0(X)$.
\end{lem}

\begin{proof}
We denote by $1$ the identity automorphism on $\wE$, and by $\theta_{\psi,\psi'} \in \mK (\mE , \mE)$, $\psi , \psi' \in \wE$, the operator defined in Rem.\ref{rem_theta}. Since $X$ is $\sigma$-compact, there is a sequence $\left\{ K_n \right\}_n$ of compact subsets covering $X$. We consider a partition of unity $\left\{ \lambda_n \right\}$ such that $\ovl{\mathrm{supp} \lambda_n} = K_n$, $n \in \bN$. By the Serre-Swan theorem, the bimodule of continuous sections of the restriction $\left. \mE \right|_{K_n}$ is finitely generated by a set $\left\{ \right. \varphi_{n,k} \left. \right\}_k$; we define $ \psi_{n,k} := \lambda_n \varphi_{n,k} \in \wE$. Furthermore, let $u_n := \sum_k \theta_{ { \psi_{n,k}} , { \psi_{n,k}} } \in \mK (\wE , \wE)$; then by construction $u_n = \lambda_n^2$ (in fact, $\sum_k \theta_{ {\varphi_{n,k}} , {\varphi_{n,k}} } = 1 |_{K_n}$) and $\sum_n^m u_n = \sum_n^m \lambda^2_n \stackrel{m}{\ra} 1$ (in the strict topology). Thus, we conclude that $\wE$ is countably generated by the set $\left\{ \psi_{n,k} \right\}$. Furthermore, if $f \in C_0 (X)$ then $\left\| f - f \sum_n^m u_n \right\| = \left\| f - \sum_n^m \lambda_n^2 f \right\| \stackrel{m}{\ra} 0$; thus $f$ is norm limit of elements of $\mK (\wE , \wE)$, and $C_0(X)$ acts on the left over $\wE$ by elements of $\mK ( \wE , \wE )$. We conclude that the pair $(\wE , 0)$ defines a class in $RK^0 (X)$.
\end{proof}

Let now $\mA$ be a $\sigma$-unital $C_0(X)$-algebra; then $\wE \otimes_X \mA$ is a countably generated $C_0(X)$-Hilbert $\mA$-bimodule. By the previous lemma, $\mA$ acts on the left on $\mM := \wE \otimes_X \mA$ by elements of $\mK (\mM , \mM)$, so that the pair $( \wE \otimes_X \mA , 0 )$ defines a class in $KK ( X ; \mA , \mA )$, corresponding to the Kasparov product $[\mE] \cdot [1]_\mA$. Thus, the 'extending the scalars' operation $\wE \mapsto \wE \otimes_X \mA$ is a particular case of the structure morphism
\begin{equation}
\label{def_i}
i_\mA : RK^0 (X) \ra  KK ( X ; \mA , \mA )
\end{equation}
\noindent defined in \cite[2.19]{Kas88}. The morphism $i_\mA$ is natural, in the following sense: if $\alpha : \mA \ra \mB$ is a $C_0(X)$-algebra isomorphism, then a ring isomorphism $\alpha_* : KK(X ; \mA , \mA) \ra KK(X ; \mB , \mB)$ is defined by the operation of passing to the pullback bimodule in the sense of Ex.\ref{ex_pullback}; thus, we find
\begin{equation}
\label{knature}
i_\mB [\mE] := [\wE {\otimes}_X \mB] 
             = [\wE {\otimes}_X \alpha (\mA)]
             = [(\wE \otimes_X \mA)_\alpha] 
             = \alpha_* i_\mA [\mE]
\end{equation}
\noindent where $(\wE \otimes_X \mA)_\alpha$ is the pullback of $\wE \otimes_X \mA$.

\bigskip

Let $\mA$ be a \sC algebra, $\mM$ a Hilbert $\mA$-bimodule; $\mM$ is said {\em imprimitivity bimodule} if there is an isomorphism $\mA \simeq \mK ( \mM , \mM )$. The set of isomorphism classes of imprimitivity bimodules, with the operation of internal tensor product, is a moltiplicative group called the {\em Picard group} of $\mA$ (see \cite[\S 3]{BGR77}: the inverse is given by the conjugate bimodule $\ovl \mM$). If $\mA$ is a $C_0(X)$-algebra, we denote by ${\bf {Pic}}_X (\mA)$ the Picard group obtained by imprimitivity $C_0(X)$-Hilbert $\mA$-bimodules. If $\mM$ is an imprimitivity $C_0(X)$-Hilbert $\mA$-bimodule, we denote by
\[  \left\{  \mM \right\} \in {\bf {Pic}}_X \mA   \]
\noindent the associated isomorphism class, and by $\left\{ \mM \right\} \cdot \left\{ \mN \right\} := \left\{ \mM \otimes_\mA \mN \right\}$ the corresponding product in the Picard group.

\begin{ex}
Let $X$ be a locally compact, paracompact Hausdorff space. The imprimitivity $C_0(X)$-bimodules correspond to modules of continuous, vanishing at infinity sections of line bundles over $X$; thus, ${\bf {Pic}}_X C_0(X)$ is isomorphic to the sheaf cohomology group $H^2 (X , \bZ)$.
\end{ex}

\bigskip

Let ${\bf {out}}_X \mA$ denote the group of $C_0(X)$-automorphisms of $\mA$ modulo inner automorphisms induced by unitaries in $M(\mA)$; let furthermore $\mK$ denote the \sC algebra of compact operators. If $\mA$ is $\sigma$-unital, then by the argument of \cite[Cor.3.5]{BGR77} we obtain a group anti-isomorphism
\begin{equation}
\label{iso_pic}
\theta : {\bf {Pic}}_X (\mA \otimes \mK) 
\stackrel{\simeq}{\longrightarrow}
{\bf {out}}_X (\mA \otimes \mK) \ \ , \ \  \left\{  \mM \right\} \mapsto \theta_\mM \ .
\end{equation}
\noindent Also recall that there is an isomorphism ${\bf {Pic}}_X (\mA) \simeq {\bf {Pic}}_X (\mA \otimes \mK)$ (see \cite[Thm.3.9]{BGR77} and successive remarks). Let $\alpha$ be a $C_0(X)$-automorphism of $\mA \otimes \mK$; then (up to inner automorphisms), it turns out that $\theta^{-1}(\alpha)$ is (the isomorphism class of) the Hilbert $(\mA \otimes \mK)$-bimodule $\mA_\alpha$ defined by
\begin{equation}
\label{def_malfa}
\mA_\alpha := \mA \otimes \mK \ \ , \ \
a \cdot \psi := \alpha (a) \psi \ \ , \ \
\psi \cdot a := \psi a \ \ , \ \ 
\left \langle \psi , \psi' \right \rangle := \psi^* \psi' \ ,
\end{equation}
\noindent where $a , \psi , \psi' \in \mA \otimes \mK$. Note that $\mO_{\mA_\alpha} \simeq \mA \rtimes_\alpha \bZ$ (\cite[\S 2]{Pim93}).

The next lemma is a version of \cite[Thm.4.2]{AEE95}:

\begin{lem}
\label{pic_pim}
Let $\mA$ be a stable $C_0(X)$-algebra, $\mM$, $\mN$ imprimitivity $C_0(X)$-Hilbert $\mA $-bimodules. Suppose that $\mM$ is {\bf outer conjugate} to $\mN$, i.e. there is $\left\{  \mV \right\} \in {\bf {Pic}}_X (\mA)$ such that $\left\{  \mM \right\} \cdot \left\{  \mV \right\} = \left\{  \mV \right\} \cdot \left\{  \mN \right\}$. Then, there is a $C_0(X)$-isomorphism $\gamma : ( \com , \bT) \ra  ( \mO_\mN  , \bT )$.
\end{lem}

\begin{proof}
Denote by $\alpha, \delta , \gamma  \in {\bf aut}(\mA)$ the automorphisms associated respectively with $\mM$, $\mN$, $\mV$, up to inner automorphisms of $\mA$. Then (up to inner automorphisms), the equality $\alpha \circ \gamma = \gamma \circ \delta$ holds. Now, by construction there are isomorphisms of Hilbert $\mA$-bimodules $\mM \simeq \mA_\alpha$, $\mN \simeq \mA_\delta$ (we use the notation (\ref{def_malfa})). Furthermore, there is a covariant isomorphism $(\beta , \gamma ) : \mA_\alpha \ra \mA_\delta$, $\beta (\psi) := \gamma (\psi)$. Thus, the lemma follows by Lemma \ref{lem_funct}.
\end{proof}

There is a natural (multiplicative) morphism
\begin{equation}
\label{def_pickk}
{\bf {Pic}}_X (\mA)  \ra  KK ( X ; \mA , \mA )  \ \ , \ \ 
\left\{  \mM \right\} \mapsto [\mM]   \ \ ,
\end{equation}
\noindent where $[\mM]$ is the class of $(\mM , 0)$ (see Rem.\ref{rem_cla}). It is clear that if $\mM$ is an imprimitivity bimodule then $[\mM] \in KK ( X ; \mA , \mA )$ is invertible; in fact, the Kasparov product reduces to the internal tensor product of Hilbert bimodules for classes of the type of Rem.\ref{rem_cla}, so that $[\mM]^{-1} = [\ovl \mM]$, $\left\{ \mM \right\} \in {\bf {Pic}}_X (\mA)$.

\section{Bimodules and endomorphisms.}
\label{bimod_alg}

The notion of {\em bimodule in a} \sC {\em algebra} has been introduced in \cite[\S 1]{DPZ97}. Let $\mB \subset \mA$ be a \sC algebra inclusion. A {\em Banach} $\mB$-{\em module in} $\mA$ is a Banach space $\mM \subset \mA$, closed for right multiplication with elements of $\mB$. If $\mM$ is also closed for left multiplication, then it is called {\em Banach} $\mB$-{\em bimodule in} $\mA$. Finally, $\mM$ is said to be {\em Hilbert} $\mB$-(bi)module in $\mA$ if
\[  
\mM^* \mM := {\mathrm {clo}}^{\left\| \cdot \right\|}
\left\{ \psi^* \psi' : \psi, \psi' \in \mM \right\} \subseteq \mB  \ , \]
\noindent (here ${\mathrm {clo}}^{\left\| \cdot \right\|}$ denotes the norm closure in $\mA$) so that $\mM$ is endowed with the $\mB$-valued scalar product $\psi , \psi' \mapsto \psi^* \psi'$. If $\mB = \mM^* \mM$, then $\mM$ is said {\em full}. Note that there is a natural identification of
\[
\mM \mM^* := {\mathrm {clo}}^{\left\| \cdot \right\|}
\left\{ \psi' \psi^* \ , \ \psi ,\psi' \in \mM \right\} \subset \mA
\]
\noindent with $\mK (\mM , \mM)$.

\bigskip

Let $\left\{ \psi_i \right\}_{i \in I}$ be a set of elements of $\mM$; we denote by $\Lambda$ the class of finite subsets of $I$. We say that $\left\{ \psi_i \right\}$ is a {\em set of generators for} $\mM$ if the net $\left\{ u_\lambda := \sum_{i \in \lambda} \psi_i \psi_i^* \right\}_{\lambda \in \Lambda}$ is an approximate unit for $\mM \mM^*$. $\mM$ is said {\em finitely generated} if there is a finite set of generators $\left\{ \psi_l \right\} \subset \mM$; in such a case, $\psi = \sum_l \psi_l (\psi_l^* \psi)$ for every $\psi \in \mM$.

A projection $P \in M(\mA)$ is called {\em support} of $\mM$ if, given an approximate unit $\left\{ u_\lambda \right\}$ for $\mM \mM^*$, the net $\left\| a u_\lambda  - aP \right\|$ converges to zero for every $a \in \mA$ (note that in such a case, the same is true for $\left\| u_\lambda a - Pa \right\|$).

It is clear that when existing, the support is unique and does not depend on the choice of the approximate unit; we will denote it by $P_\mM$. If $\mM$ is finitely generated by a set $\left\{ \psi_l \right\}$ then $P_\mM = \sum_l \psi_l \psi_l^*$, so that $P_\mM \in \mA$.

\bigskip

\begin{rem}
Our definition of support is different from the one given in \cite[\S 2]{CDPR94}, where just the convergence in the strong topology is required. An example is the Hilbert space of isometries $H := {\mathrm{span}}\left\{ \psi_k \right\}_{k \in \bN}$ generating the Cuntz algebra $\mO_\infty$ (regarded as a Hilbert $\bC$-bimodule in $\mO_\infty$): note in fact that the sequence $\sum_k^n \psi_k \psi_k^*$ does not converge in the strict topology. If $\mM$ is finitely generated, our notion of support coicide with the one given in \cite[\S 2]{CDPR94}.
\end{rem}

\begin{ex}
\label{ex_sk}
Let $H$ be a separable Hilbert space, $B(H)$ the \sC algebra of bounded linear operators on $H$; it is well known that $B(H) = M(\mK)$. We regard $\mK$ as a Hilbert $\mK$-bimodule in $B(H)$. If $\left\{ e_k \right\}_{k \in \bN}$ is an orthonormal basis for $H$, we consider the operators $\theta_k := \theta_{e_k , e_k} \in \mK$, with the notation Rem.\ref{rem_theta}. Note that $\theta_k = \theta_k^* = \theta_k^2$. Now, the sequence $\left\{ u_n = \sum_{k=1}^n \theta_k \theta_k^* \right\}_n$ converges to the identity in the strict topology, so that $\mK$ has support $1$ in $B(H)$.
\end{ex}

\begin{ex}
\label{ex_son}
We consider the Cuntz algebra $\mO_d$, $d \in \bN$. The Banach space $\mO_d^1$ is a Hilbert $\mO_d^0$-bimodule in $\mO_d$ (we use the notation of \S \ref{pimsner}). Since $\mO_d^1$ contains the generators $\psi_1 , \ldots , \psi_d$ of $\mO_d$, we find that $\mO_d^1$ has support $1 = \sum_i^d \psi_i \psi_i^* $. Let us now consider the stable Cuntz algebra $\mO_d \otimes \mK$; then, $\mO_d^1 \otimes \mK$ is a Hilbert $(\mO_d^0 \otimes \mK)$-bimodule in $\mO_d \otimes \mK$; with the notation of the previous example, we find that $\left\{ \psi_i \otimes \theta_k \right\}_{i,k}$ is a set of generators for $\mO_d^1 \otimes \mK$; thus, $\mO_d^1 \otimes \mK$ has support $1$.
\end{ex}

The following lemma is a version of \cite[Thm.3.12]{Pim93}.

\begin{lem}
\label{lem_bima}
Let $\mA$ be a \sC algebra, $\mM$ a full Hilbert $\mB$-bimodule in $\mA$ with support $1$. Then, the inclusion $\mM \subset \mA$ defines a \sC algebra morphism $i_\mM : \com \ra \mA$, such that
\[  i_\mM (b \psi b') = b \ i_\mM (\psi) \ b'  \ \ , \ \ 
    i_\mM \left \langle \psi , \psi' \right \rangle = i_\mM (\psi)^* i_\mM (\psi') \]
\noindent where $b , b' \in \mB$, $\psi , \psi' \in \mM$, and $\left \langle \cdot , \cdot \right \rangle$ is the $\mB$-valued scalar product. 
\end{lem}

\begin{rem}
The analogue of the previous lemma holds for the Toeplitz algebra $\mT_\mM$ (see \cite{Pim93} for the definition of $\mT_\mM$, or Sec.\ref{basics} of the present paper) instead of $\com$, if $\mM$ does not have support the identity (see \cite[Thm.3.4]{Pim93}).
\end{rem}

\bigskip

Let $\mZ$ denote the centre of $\mA$. If $\rho$ is an endomorphism on $\mA$, we consider for every $r,s \in \bN$ the intertwiners spaces 
\[ 
(\rho^r , \rho^s) := \left\{ t \in \mA : t \rho^r(a) = \rho^s (a) t , a \in \mA \right\} \ .
\]
\noindent Every $(\rho^r , \rho^s)$ has a natural structure of Banach $\mZ$-bimodule in $\mA$. In particular, if we consider the identity automorphism $\iota$, then we obtain that every $(\iota , \rho^r)$ has a further structure of Hilbert $\mZ$-bimodule, by defining the $\mZ$-valued scalar product $t , t' \mapsto t^* t'$ for $t , t' \in (\iota , \rho^r)$. Let us now consider the multiplier algebra $M(\mA)$, with the centre $ZM(\mA)$. We consider the Banach space
\begin{equation}
\label{def_innmod}
M_\rho := \left\{ \psi \in M(\mA) : \psi a = \rho (a) \psi , a \in \mA \right\}   \ ,
\end{equation}
\noindent and the abelian \sC algebra
\[
\zro := {\mathrm {clo}}^{\left\| \cdot \right\|}
\left\{ \psi^* \psi' \ : \  \psi , \psi' \in M_\rho \right\} \subseteq ZM(\mA) 
\]
\noindent (note that $\psi^* \psi' \in \mA' \cap M(\mA) = ZM(\mA)$, $\psi , \psi' \in M_\rho$). By construction, $M_\rho$ is a full Hilbert $\zro$-bimodule in $M(\mA)$. 

\begin{defn}
\label{inner}
Let $\rho$ be an endomorphism of a \sC algebra $\mA$. $\rho$ is said to be {\bf topologically inner} if $M_\rho$ has support $P_\rho$, and $\rho (a) P_\rho = \rho (a)$ for every $a \in \mA$. If $M_\rho$ is finitely generated, then $\rho$ is said {\bf inner}.
\end{defn}

\noindent When $\rho$ is topologically inner, then for every set $\left\{  \psi_i \right\}$ of generators for $M_\rho$ the following equality holds:
\begin{equation}
\label{rel_inn}
\rho (a) = \rho (a) P_\rho = \rho (a) \ {\mathrm {Lim}}_\lambda u_\lambda =
           {\mathrm {Lim}}_\lambda 
           \sum_{i \in \lambda} \psi_i a \psi_i^* \ \ , \ a \in \mA \ . 
\end{equation}
\noindent Note that in the case in which $\mA$ is unital, then $M_\rho = (\iota , \rho)$.

\begin{rem}
Let $v \in M(\mA)$ be a partial isometry, $v^*v = 1$, $vv^* = P$; then, an inner endomorphism $\rho (a) := v a v^*$, $a \in \mA$ is induced, and $M_\rho$ is a free $ZM(\mA)$-bimodule in $M(\mA)$ with support $P$; thus, the usual notion of inner endomorphism is recovered. 
\end{rem}

\begin{prop}
\label{rem_proj_end}
The \sC algebra $\zro$ is a (closed) ideal of $ZM(\mA)$. If $\rho$ is inner, then $M_\rho$ is a finitely generated, projective Hilbert $ZM(\mA)$-bimodule.
\end{prop}

\begin{proof}
Since $M_\rho$ is closed for multiplication by elements of $ZM(\mA)$, we obtain that $\psi^* \psi' z \in \zro$ for $\psi , \psi' \in M_\rho$, $z \in ZM(\mA)$. Thus $\zro$ is an ideal of $ZM(\mA)$. 

Let now $\rho$ be inner. Then, $M_\rho$ is finitely generated by a set $\left\{ \psi_l \right\}_{l=1}^n$. In order to check the projectivity, we consider the matrix $E := (\psi_l^* \psi_m)_{l,m} \in ZM(\mA) \otimes \bM_n$ and note that
\[
(E^*)_{lm} = (E_{ml})^* = \psi_m^* \psi_l = E_{lm} \ \ , \ \ 
(E^2)_{lm} = \sum_k \psi_l^* \psi_k \psi_k^* \psi_m = \psi_l^* P_\rho \psi_m = E_{lm} .
\]
\noindent Thus, $E$ is a projection in $ZM(\mA) \otimes \bM_n$. Let $M_E := E (\oplus^n ZM(\mA))$ be the projective Hilbert $ZM(\mA)$-module associated with $E$. Then, the map $\beta : M_\rho \ra M_E$, $\beta (\psi) := (\psi_l^* \psi)_l$ supplies the desired unitary isomorphism.
\end{proof}

\begin{rem}
Let $\mA$ be a \sC algebra, $\mB \subseteq ZM(\mA)$, $\mM \subset M(\mA)$ a full Hilbert $\mB$-bimodule with support $P_\mM$. If $\left\{ \psi_\alpha \right\}$ is a set of generators for $\mM$, then a topologically inner endomorphism $\rho$ is induced on $\mA$, defined as by (\ref{rel_inn}). It is clear that $\rho$ does not depend on the choice of $\left\{ \psi_\alpha \right\}$. By construction, $M_\rho = \mM \cdot ZM (\mA) \simeq \mM \otimes_\mB ZM (\mA)$.
\end{rem}

\begin{rem}
Let $X$ be the spectrum of $ZM(\mA)$. Suppose $P_\rho = 1$; then $M_\rho \cdot M_\rho^*$ is a $C(X)$-operator system in the sense of \cite[Appendix]{Bla97}.
\end{rem}

\begin{ex}
The canonical endomorphism $\sigma$ of the Cuntz algebra $\mO_d$, $d < \infty$, is inner in the sense of the previous definition (see \cite{DR87,Cun81}); $M_\sigma$ coincides in fact with the Hilbert space of isometries generating $\mO_d$.
\end{ex}

\begin{ex}
Let $\mA$ be a continuous trace \sC algebra, $\rho$ a locally unitary automorphism in the sense of \cite{PR84}. Then $\rho$ is topologically inner in the sense of Def.\ref{inner}; $M_\rho$ coincides with the module of continuous sections of the line bundle associated with $\rho$.
\end{ex}

\begin{ex}
We generalize the previous example. Let $X$ be a paracompact Hausdorff space, $\mA$ a stable continuous trace \sC algebra with spectrum $X$. Then, every $C(X)$-endomorphism of $\mA$ is topologically inner in the sense of Def.\ref{inner} (see \cite{Hir01}). If $X$ is compact, then $\rho$ is inner.
\end{ex}

We will see in the sequel that a fundamental class of examples for Def.\ref{inner} is given by the canonical endomorphisms over Pimsner (Toeplitz) algebras of vector bundles.

Crossed product constructions, performed in such a way that an endomorphism becomes inner in the sense of Def.\ref{inner}, are made in \cite[\S 3]{Vas03}.

\section{Basic Properties of $\coe$.}
\label{basics}

Let $X$ be a locally compact Hausdorff space, $\mH \ra X$ a rank $d$ Hilbert bundle. For every $x \in X$, we denote by $\mH_x \simeq \bC^{d(x)}$ the fibre of $\mH$ over $x$. In the present section we give some basic properties of the Pimsner (DR) algebra associated with the Hilbert $C_0(X)$-bimodule $\wH$ of continuous, vanishing at infinity sections of $\mE$. The left $C_0(X)$-action of $\wH$ is defined to coincide with the right one.

We denote by $\coh$ the Pimsner algebra of $\wH$, and by $\doh$ the DR-algebra of $\wH$. In particular, the analogous notations will be used for a vector bundle $\mE \ra X$.

Let $(\hrs)$ denote the set of bounded Hilbert bundle morphisms from the tensor power $\mH^r$ into $\mH^s$, $r,s \in \bN$. Then, $(\hrs)$ is a Banach $C_0(X)$-bimodule, which by definition (\cite[I.4]{DD63}) coincides with the space $(\wH^r , \wH^s)$ of Hilbert $C_0(X)$-bimodule operators between the corresponding tensor powers of $\wH$. In the sequel, we will make the identification 
\[
(\wH^r , \wH^s) \equiv (\hrs) \ \ ,
\]
\noindent so that our terminology for $\coh$ (resp. $\doh$) as the Pimsner (DR) algebra of the Hilbert bundle $\mH$ is justified. It follows from the previous identity that, in particular, the spaces of compact operators $\mK (\wH^r , \wH^s)$ are characterized in terms of Hilbert bundle morphisms. According to such identification, we will use the notation $\mK (\hrs) := \mK (\wH^r , \wH^s)$. We also adopt the notation $\iota := \mH^0 := X \times \bC$, so that $\wH \simeq \mK (\iota , \mH)$. Note that $(\iota , \iota) = C_b (X)$, $\mK (\iota , \iota) = C_0 (X)$.

\bigskip

\begin{rem}
\label{h_fin_gen}
When $X$ is compact and $\mE \ra X$ is a vector bundle it is well-known that $\wE$ is finitely generated, so that $\mK (\ers) = (\ers)$ for every $r,s \in \bN$. Thus, in such a case $\doe$ coincides with $\coe$.
\end{rem}

\bigskip

Let us now discuss the natural immersion $\wH \hra (\iota,\mH)$. If $\psi : X \ra \mH$ is a continuous, vanishing at infinity section belonging to $\wH$, we can define a Hilbert bundle morphism $\widetilde \psi$ from $X \times \bC$ into $\mH$, by assigning $\widetilde \psi (x , z) := z \cdot \psi (x)$. Viceversa, every Hilbert bundle morphism $\varphi : X \times \bC \ra \mH$ defines a continuous section $\varphi_0 : \varphi_0 (x) := \varphi (x , 1)$. Anyway, $\varphi_0$ does not necessarily vanish at infinity: in order to have a well defined Hilbert bundle morphism, it suffices in fact that $\sup_{x} \left\| \varphi (x , 1) \right\| < \infty$. Note that $\varphi$ induces a Hilbert $C_0(X)$-bimodule morphism from $\wH$ onto $C_0(X)$, in the following way: if $\psi \in \wH$, we denote by $\left \langle \varphi_0 , \psi \right \rangle$ the continuous function $\left\{ x \mapsto  {\left \langle \varphi_0 (x) , \psi (x) \right \rangle}_x \right\}$, where ${\left \langle \cdot , \cdot \right \rangle}_x$ denotes the $\bC$-valued scalar product on $\mH_x \simeq \bC^{d(x)}$. Then $\left| {\left \langle \varphi_0 (x) , \psi (x) \right \rangle}_x \right| \leq \left\| \varphi_0 (x) \right\| \left\| \psi (x) \right\|$, so that $\left \langle \varphi_0 , \psi \right \rangle$ actually vanishes at infinity (since $\left\| \psi (x) \right\|$ vanishes at infinity). Note that in the setting of the DR-algebra $\doh$, we have $\varphi_0 \in (\iota , \mH)$, and the map $\left\{ \psi \mapsto {\left \langle \varphi_0 , \psi \right \rangle} \right\}$ corresponds to $\varphi_0^* \in (\mH , \iota)$.

\bigskip

\begin{rem}
\label{full_bimod}
In the sequel, we will make use of the notion of total subset of a continuous field of Banach spaces. Let $\mH \ra X$ be a Banach bundle; if $M \subset \wH$ is a total subset for the continuous field of Banach spaces associated with $\mH$, then the closed $C_0(X)$-bimodule generated by $M$ coincides with $\wH$ (see \cite[10.2.3]{Dix}).
\end{rem}

\bigskip

Let now $\mH^* \ra X$ be the dual Hilbert bundle. We consider the Banach bundle
\[  p_{rs} : \mH^{r,s} := \mH^s \otimes {\mH^*}^r \ra X  \ , \]
\noindent having fibre the space $\bM_{ d(x)^r , d(x)^s }$ of $d(x)^r \times d(x)^s$ matrices. Now, with the notation of Rem.\ref{rem_theta}, we obviously find 
\[  
\mK (\hrs) := \left\{ \theta_{\psi , \psi'} , \psi \in \wH^s , \psi' \in \wH^r \right\} \subseteq C_0(X , \mH^{r,s}) \ \ .
\]
\noindent If $w \in \mH^{r,s}$, $p_{rs}(w) = x$, then we can write (up to linear combinations) $w = \theta_{v , v'}$, $v \in \mH^s_x \simeq \bC^{d(x)^s}$, $v' \in \mH^r_x \simeq \bC^{d(x)^r}$; since $\mH^r$, $r \in \bN$, are full Hilbert bundles we conclude that $w = \theta_{\psi,\psi'} (x)$, where $\psi \in \wH^s$, $\psi' \in \wH^r$ and $\psi (x) = v$, $\psi'(x) = v'$. Thus, for every $w \in \mH^{r,s}$ there is $t \in \mK (\hrs)$ such that $t (x) = w$; this property implies that $\mH^{r,s}$ defines the unique continuous field of Banach spaces admitting $\mK (\hrs)$ as a total subset (see \cite[10.2.3]{Dix}). So that, we have the following

\begin{lem}
\label{K-bundle}
Every $\mK (\hrs)$, $r,s \in \bN$, coincides with the module $C_0 (X, \mH^{r,s})$ of continuous, vanishing at infinity sections of the Banach bundle $\mH^{r,s} \ra X$. In particular, $\mK (\mH^r , \mH^r)$ is a continuous bundle of \sC algebras, with fibres the order $d(x)^r$ matrix algebras $\bM_{d(x)^r}$.
\end{lem}

\begin{ex}
\label{trivial}
Let $X$ be a compact Hausdorff space, $\mE := X \times \bC^d$ the trivial vector bundle with rank $d > 1$. Then $\doe = \coe \simeq C(X) \otimes \mO_d$, where $\mO_d$ is the order $d$ Cuntz algebra. In fact, $\wE$ is the free rank $d$ Hilbert $C(X)$-bimodule, so that $(\ers) \simeq C(X) \otimes \bM_{d^r , d^s}$ for every $r,s \in \bN$. In the case in which $\mE$ has rank $1$, recall that the Pimsner algebra associated with the free rank $1$ Hilbert bimodule $\bC$ is given by $C(S^1)$, that we regard as the \sC algebra generated by a unitary $u$ such that $u^r \neq 1$ for every $r \in \bN$. So that, the isomorphism $\mE \simeq X \times \bC$ gives $\coe \simeq C(X) \otimes C(S^1)$, with the same argument for $d > 1$.
\end{ex}

The following result has been proved in the abstract setting of certain tensor \sC categories in \cite[Thm.4.1]{Vas} for 'locally trivial objects', but previously (and independently) by J.E. Roberts (\cite{Rob}) in the case of vector bundles.

\begin{prop}
\label{coe_cf}
Let $X$ be a locally compact Hausdorff space, $\mH \ra X$ a Hilbert bundle. Then $\coh$ is a continuous bundle of Cuntz algebras over $X$, and $\doh$ defines a continuous field of Cuntz algebras over $X$. If $\mE \ra X$ is a rank $d$ vector bundle, then $\coe$ is an $\mO_d$-bundle (i.e., $\coe$ is locally trivial). If $\mE$ is a line bundle, then $\coe$ is a $C(S^1)$-bundle, where $C(S^1)$ is the \sC algebra of continuous functions over the circle.
\end{prop}

\begin{proof}
Let $d(x)$ denote the rank of $\mH$ over $x \in X$. Then, every $t \in (\hrs)$ defines a vector field $\left\{ t_x \in \bM_{d(x)^r  , d(x)^s} \subset \mO_{d(x)} : x \in X \right\}$, and the norm function $\left\{ x \mapsto \left\| t_x \right\| \right\}$ is bounded and continuous. We define  
\[
\widetilde \Theta_\mH := \left\{ \left\{ t_x \right\} \in \prod_x \mO_{d(x)} : t \in \ ^0 \doh \right\} \ .
\] 
\noindent By the previous lemma, we find that for every $x \in X$, $t_0 \in \ ^0 \mO_{d(x)}$ there is $t \in \widetilde \Theta_\mH$ such that $t_0 = t_x$. Thus, by applying \cite[10.2.3]{Dix}, we find that  $\widetilde \Theta_\mH$ defines a continuous field of \sC algebras. Now, for every $t \in  \Theta_\mH$ the $\sup${\em -norm} $\sup_x \left\| t_x \right\|$ is well defined, and the circle action (\ref{def_circ_act}) is isometric w.r.t. such a norm. By universality of the DR-algebra, we find that $\doh$ can be identified as the \sC algebra obtained by the closure of $\widetilde \Theta_\mH$ w.r.t. the $\sup$-norm.

\noindent We emphasize the fact that in general $\coh$ does not coincide with the set of continuous, vanishing at infinity vector fields of $\widetilde \Theta_\mH$. We will return on this point in Rem.\ref{rem_hbundle}.

\noindent We now pass to the Pimsner algebra $\coh$. First note that $\coh$ is obviously a Banach $C_0(X)$-bimodule for left and right multiplication. Furthermore,  by Lemma \ref{K-bundle}, if $t \in \mK (\hrs)$ then the norm function vanishes at infinity. With the same argument used for the DR-algebra, we find that $\coh$ can be regarded as the \sC algebra of the continuous field $\Theta_\mH := \left\{ \left\{ t_x \right\} \in \prod_x \mO_{d(x)} : t \in \ ^0 \coh \right\}$. Thus, $\coh$ is a continuous bundle of Cuntz algebras.

\noindent Let now $U \subset X$ be a closed set trivializing $\mH$; for a generic (non-locally trivial) Hilbert bundle, we can pick $U = \left\{ x \right\}$, $x \in X$. We consider the corresponding local chart $\pi_U : \left. \mH \right|_{U} \ra U \times \bC^{d(x)}$, $x \in U$. Then, a covariant morphism $\pi_U^* : \wH \ra C(U) \otimes \bC^{d(x)}$ is induced, with $\pi_U^* (f \psi) = \left. f \right|_{U} \pi_U^*(\psi)$, $f \in C_0(X)$, $\psi \in \wH$. It follows from Lemma \ref{lem_funct} and Ex.\ref{trivial} that a local chart
\[  {\wa \pi}_U :  \doh |_{U} \ra C(U) \otimes \mO_{d(x)}  \]
\noindent is induced. The existence of such local charts shows that the fibre of $\coh$ over $x$ is the Cuntz algebra $\mO_{d(x)}$, and that $\coh$ is locally trivial if $\mH$ is a vector bundle. The proof in the case $d(x) \equiv 1$ is analogue, with $C(S^1)$ playing the role of the Cuntz algebra.
\end{proof}

\begin{rem}
\label{rem_line_bundle}
Recall that if $\mE \ra X$ is a vector bundle, the {\em sphere bundle} of $\mE$ is the topological bundle $\mE_1 \ra X$, having as fibre the space of elements of $\mE$ with norm $1$. In particular, if $\mL \ra X$ is a line bundle, then $\mL_1 \ra X$ has fibre $S^1$. It turns out that $\mO_\mL$ is abelian, with spectrum given by the sphere bundle $\mL_1$. Furthermore, if $\mL^* \ra X$ is the dual bundle of $\mL$, then $\mO_{\mL^*}$ is isomorphic to $\mO_\mL$ (see \cite[Prop.4.3]{Vas}). Note that in general $\mO_\mL \simeq C_0(\mL_1)$ is a nontrivial $C(S^1)$-bundle (in fact, in general $\mL_1 \ra X$ is nontrivial as a topological bundle).
\end{rem}

\begin{rem}
\label{nontrivial}
If the vector bundle $\mE \ra X$ has rank $d > 1$, the fact that $\coe$ in general is a {\em nontrival} $\mO_d$-bundle can be checked by calculating the $K$-theory of $\coe$, as done in \cite{Gol00} in the case in which $X$ is a $n$-sphere, and in \cite{Vasu} in the case in which $X$ is a Riemannn surface. These results will be exposed in a future paper.
\end{rem}

\begin{rem}
\label{nuclear}
The previous proposition implies that $\coh$ is nuclear. In fact, $\coh$ is a continuous bundle having as fibres nuclear \sC algebras (see \cite[\S 2.2]{Bla95}). In particular, $\coh$ is also exact (as can be concluded also by \cite[Cor.2.2]{DS99})
\end{rem}

We denote by $\wa \mO_\mH \ra X$ the {\em topological bundle associated with} $\coh$, as defined in Sec.\ref{fields}. Note that for every $r,s \in \bN$, there is a bundle inclusion $\mH^{r,s} \hra \wa \mO_\mH$. The following result is an immediate consequence of \cite[Thm.3.3]{APT73}.

\begin{cor}
Let $\mE \ra X$ be a vector bundle. Then, the multiplier algebra $M(\coe)$ is isomorphic to the \sC algebra $C_b (X , \wa \mO_\mE )$ of bounded continuous sections of $\wa \mO_\mE$.
\end{cor}

\begin{rem}
In the case of a Hilbert bundle, \cite[Thm.3.3]{APT73} again implies that $M(\coh)$ is the \sC algebra of bounded continuous sections of a topological bundle of Cuntz algebras over $X$. Anyway, the topology of such a bundle is different w.r.t. the one of $\wa \mO_\mH \ra X$ (see Rem.\ref{rem_hbundle} below).
\end{rem}

\bigskip

Let $\mE \ra X$ be a rank $d$ vector bundle. The action introduced in Cor.\ref{u_act} defines a canonical $\ud$-action by automorphisms over the Cuntz algebra $\mO_d$. Such an action can be used to determine the ${\bf aut}\mO_d$-cocycle associated with $\coe$ as an $\mO_d$-bundle, in the case in which $X$ is paracompact (see \cite[\S 4]{Vas} for details). Let in fact $u := ( \left\{ X_i\right\} ,\left\{ u_{ij} \right\} ) \in H^1 (X , \ud)$ be an $\ud$-cocycle associated with $\mE$ (i.e., $\left\{ X_i\right\}$ is a locally finite open cover, and every $u_{ij}$ is a continuous map from $X_i \cap X_j$ into $\ud$, satisfying the relation $u_{ik} = u_{ij} u_{jk}$; we use notation and terminology of \cite{Kar}). As we noted in the proof of the previous proposition, local charts $\pi_i : \left. \mE \right|_{X_i} \ra X_i \times \bC^d$ induce by functoriality local charts ${\wa \pi}_i : \left. \coe \right|_{X_i} \ra C_0 (X_i) \otimes \mO_d$. Since up to cocycle equivalence $u_{ij} (x) = \pi_i (x) \cdot \pi_j^{-1} (x) \in \ud$, again by functoriality we find the equality $\wa u_{ij} (x) = {\wa \pi}_i (x) \circ {\wa \pi}_j^{-1} (x)$, where $\wa u_{ij} (x) \in {\bf aut} \mO_d$ is the automorphism defined as by Cor.\ref{u_act}. The ${\bf aut}\mO_d$-cocycle associated with $\coe$ is then given by $\wa u := ( \left\{ X_i \right\} ,\left\{ \wa u_{ij} \right\} ) \in H^1 (X , {\bf aut} \mO_d )$. It follows from the previous considerations that the bundle $\wa \mO_\mE \ra X$ can be easily  described by clutching the topological spaces $X_i \times \mO_d$, {\em via} the transition maps $\wa u_{ij} : X_i \cap X_j \ra {\bf aut} \mO_d$.

\begin{ex}
\label{ex_line_bun}
Let $X$ be a paracompact Hausdorff space, $d \in \bN$, $d > 1$. Since $\bT$ is a subgroup of $\ud$, there is a map $i : H^1 ( X , \bT  ) \ra H^1 ( X , \ud  )$. It follows from the above considerations that every cocycle $z \in H^1(X ,\bT)$ defines an $\mO_d$-bundle $\mO_{z,d}$ over $X$, with cocycle $\wa {i(z)} \in H^1 (X, {\bf aut} \mO_d)$. Now, $z$ can be regarded as a set of transition maps $z_{ij} : X_i \cap X_j \ra {\bT}$ for a line bundle $\mL \ra X$; in such a way, $i(z) \in H^1 (X , \ud)$ is an $\ud$-cocycle for the vector bundle $\eps (d) \otimes \mL$, where $\eps (d) := X \times \bC^d$. Thus $\mO_{z,d}$ is the Pimsner algebra associated with $\eps (d) \otimes \mL$. In other terms, we defined a map from $H^1 (X , \bT)$ into the set of isomorphism classes of $\mO_d$-bundles over $X$.
\end{ex}

\begin{prop}
\label{coe_doe}
Let $X$ be a locally compact Hausdorff space, $\mE \ra X$ a vector bundle. Then, there are \sC algebra isomorphisms 
\[  
\doe \simeq  M(\coe) \simeq C_b ( X , \wa \mO_\mE ) \ .
\]
\end{prop}

\begin{proof}
It follows from the elementary theory of vector bundles that $(\ers) \simeq C_b (X , \mE^{r,s})$ (see \cite[I.5.9]{Kar}). Since there is a bundle inclusion $\mE^{r,s} \hra \wa \mO_\mE$, every element $t$ of $^0\doe$ defines a vector field $\left\{ t_x \in \ ^0\mO_d \right\}_x$ belonging to $C_b (X , \wa \mO_\mE )$. Let us now consider the $\sup$-norm over $C_b (X , \wa \mO_\mE )$; it is clear that the circle action is isometric w.r.t. such norm, so that by universality there is a \sC algebra monomorphism $\doe \hra C_b ( X , \wa \mO_\mE )$. Now, the \sC algebra $C_b ( X , \wa \mO_\mE )$ defines a continuous field of \sC algebras, with fibre $\mO_d$: since for every $x \in X$ the set $\left\{ t_x , t \in \doe \right\}$ is dense in $\mO_d$, we find that $\doe$ coincides with $C_b (X , \wa \mO_\mE)$, and the proposition is proved.
\end{proof}

Thus $\coe$ and $\doe$ describe essentially the same object, the topological bundle $\wa \mO_\mE \ra X$. It is clear that every \sC algebra bundle property of $\coe$ can be extended to $\doe$. For example, an isomorphism $\coe \simeq \mO_{\mE'}$ of continuous bundles is equivalent to an isomorphism $\doe \simeq \widetilde \mO_{\mE'}$.

\bigskip

\begin{rem}
\label{rem_hbundle}
The previous proposition fails for a generic Hilbert bundle $\mH \ra X$. In order to prove this assertion, we expose a class of examples (see \cite[\S 5]{Dup74}). Let $X$ be a compact Hausdorff space, $A \subset X$ a closed subset; if $d \in \bN$, we consider the Hilbert bundle $\mH := ((X - A) \times \bC^{d} ) \sqcup (A \times \bC^{d+1})$. The corresponding module of continuous sections is given by $\wH := C^d(X) \oplus C_A (X)$, where $C^d (X) := \left\{ f : X \ra \bC^d : f \ {\mathrm {continuous}} \right\}$, $C_A (X) := \left\{ f \in C(X) : \left. f \right|_{X - A} = 0 \right\}$ are endowed with the obvious Hilbert $C(X)$-bimodule structures. The \sC algebra $\mK ( \mH , \mH )$ is a non-unital continuous bundle over $X$; in fact, we can describe $\mK ( \mH , \mH )$ in terms of the matrix algebra
\[
\mK ( \mH , \mH ) \simeq
\left(
\begin{array}{cc}
    C(X,\bM_d) & 
    \mK (C^d(X) , C_A(X)) 
 \\ \mK (C_A(X) , C^d(X))    & 
    C_A(X)
\end{array}
\right)
\]
%
%
\noindent that does not contain an identity. Let now $\chi_A$ be the characteristic function of $A$. Then, the matrix $I := {\mathrm {diag}}(1_d , \chi_A)$ is a multiplier of $\mK ( \mH , \mH )$ (in fact, $Ia = aI = a$, $a \in \mK ( \mH , \mH )$), and belongs to $( \mH , \mH )$ (note that $I \psi = \psi$ for every $\psi \in \wH$). Thus, $I$ is an element of $(\mH , \mH)$, but not of $\mK ( \mH , \mH )$. In other terms, $I \in \doh$ but is not a continuous section of $\wa \mO_\mH$.
\end{rem}

\begin{lem}
\label{generators}
Let $X$ be a locally compact, paracompact Hausdorff space, $\mE \ra X$ a vector bundle with rank $d$. Then, $\wE$ is a Hilbert $C_0(X)$-bimodule in $\coe$ having as support the identity, i.e. $\wE$ admits a set of generators in the sense of Sec.\ref{bimod_alg}. The same is true for every $\wE^r$, $r \in \bN$.
\end{lem}

\begin{proof}
Let us consider 

\begin{itemize}

\item  a locally finite, open trivializing cover $\mU := \left\{ U_i \right\}_{i \in I}$ for $\mE$;

\item  a partition of unity $\left\{ \chi_i \right\}$, subordinate to $\mU$ (so that, $\sum_i \chi_i^2$ converges to $1$ in the strict topology);

\item a set $\left\{ \psi_{i,k} \in \wE \ , \ k = 1, \cdots , d \ , \ i \in I \right\}$, such that $\psi_{i,h}^* \psi_{i,k} = \delta_{hk} \chi^2_i$, $h,k = 1 , \cdots , d$ (here $\delta_{hk}$ denotes the Kronecker symbol), and $\sum_k \psi_{i,h} \psi_{i,h}^* = \chi^2_i$. Such a set exists by local triviality.

\end{itemize}

\noindent Let $\Lambda$ denotes the class of finite subsets of $I$. We consider the net
\[  
\left\{ u_\lambda := 
\sum \limits_{i \in \lambda} \sum \limits_k^d \psi_{i,k} \psi_{i,k}^* 
\right\}_{\lambda \in \Lambda} \ \ ;
\]
\noindent note that if $\psi \in \wE$, then $\sum_{k} \psi_{i,k} \psi_{i,k}^*  \psi = \chi_i^2 \psi$. Thus, $\left\{ u_\lambda \right\}$ converges in the strict topology to the identity of $(\mE , \mE)$, and is an approximate unit for $\coe$. The assertion for $\wE^r$ follows by considering sets of generators of the type $\left\{ \psi_{i,k_1} \cdots \psi_{i,k_r}  \right\} \subset (\iota , \mE^r)$, after rescaling of the partitions of unity.
\end{proof}

\bigskip

Let $\coh^k$, $k \in \bZ$, denote the spectral subspaces w.r.t. the circle action (\ref{def_circ_act}). From Prop.\ref{coe_cf}, we have the following

\begin{cor}
For every $k \in \bZ$, the spectral subspace $\coh^k$ is a Banach bundle, with fibres $\mO_{d(x)}^k$. In particular, $\coh^0$ is a \sC algebra bundle with fibre the {\it {UHF}} algebras $\mO_{d(x)}^0$. All the bundles are locally trivial if $\mH$ is a vector bundle. In the same way, every $\doh^k$ is a continuous field of Banach spaces.
\end{cor}

\begin{rem}
\label{hilb_bundle}
Every space $\coh^k$ has the following natural structure of a $C_0(X)$-Hilbert $\coh^0$-bimodule:
\[  
a,t \mapsto at  \ \ , \ \  
t,a \mapsto ta \ \ , \ \ 
\left \langle  t , t' \right \rangle := t^* t'
\]
\noindent where $a \in \coh^0$, $t , t' \in \coh^k$. In particular, the Hilbert $\coh^0$-bimodule $\coh^1$ has been considered in \cite[\S 2]{Pim93} (denoted there by $E_\infty$, where $E$ is the Hilbert bimodule generating the Pimsner algebra). Let $k \geq 0$; then, the natural map $\beta : \wH^k \otimes_X \coh^0 \ra \coh^k$, $\beta (\psi \otimes a) := \psi a$, defines an isomorphism of right Hilbert $\coh^0$-bimodules. Of course, the same considerations apply to the $C_b(X)$-Hilbert $\doh^0$-bimodules $\doh^k$.
\end{rem}

\bigskip

\begin{rem}
\label{zero_grade}
Let $\mE \ra X$ be a vector bundle. Then, for every $r \in \bN$ there is a $C_0(X)$-algebra isomorphism $\mK (\mE^r ,\mE^r) \simeq \otimes_X^r \mK (\mE , \mE)$, so that the zero-grade bundle $\coe^0$ has a filtration 
\begin{equation}
\label{filtration}
\coe^0 \simeq \lim \limits_{\ra r} \otimes_X^r \mK (\mE , \mE) \ \ .
\end{equation}
\noindent Let $\mE' \ra X$ be a vector bundle such that $\mE^r$ is isomorphic to $\mE'^r$ for some $r \in \bN$. Then, for every line bundle $\mL \ra X$, there are isomorphisms $\coe^0 \simeq \mO_{\mE' \otimes \mL}^0 \simeq \mO_{\mE'}^0$. In fact, there are isomorphisms of matrix algebra bundles 
\[
\mK ( (\mE' \otimes \mL)^r , (\mE' \otimes \mL)^r ) 
\stackrel{\alpha_r}{\longrightarrow} 
\mK ( \mE'^r , \mE'^r )
\stackrel{\alpha'_r}{\longrightarrow} 
\mK ( \mE^r , \mE^r ) \ ;
\]
\noindent that we extend to isomorphisms 
$\mO_{\mE' \otimes \mL}^0
\stackrel{\alpha}{\longrightarrow} 
\mO_{\mE'}^0
\stackrel{\alpha'}{\longrightarrow} 
\coe^0$ by defining (thanks to the filtration (\ref{filtration}))
\[
\begin{array}{ll}
\alpha 
    (t_{11} \otimes \cdots \otimes t_{1r} \otimes 
     t_{21} \otimes \cdots \otimes t_{2r}  \otimes \cdots ) := \\
\alpha_r (t_{11} \otimes \cdots \otimes t_{1r}) \otimes 
\alpha_r (t_{21} \otimes \cdots \otimes t_{2r}) \otimes  \cdots \ \ ,
\end{array}
\]
\noindent $t_{ij} \in \mK (\mE' \otimes \mL , \mE' \otimes \mL)$. The isomorphism $\alpha'$ is constructed in the same way. In particular, for every vector bundle $\mE \ra X$ and line bundle $\mL \ra X$, there is an isomorphism $\coe^0 \simeq \mO_{\mE \otimes \mL}^0$.
\end{rem}

\bigskip

We now give some elementary properties of the \sC algebras $\coh$ w.r.t. continuous maps over the base spaces. We recall that if $p : Y \ra X$ is a continuous surjective map of locally compact Hausdorff spaces, then for every Hilbert bundle $\mH \ra X$ the {\em pullback bundle} $p_* \mH \ra Y$ is given by
\[  p_* \mH := \mH \times_X Y := \left\{ (v , y) \in \mH \times Y : v \in \mH_{p(y)} \right\}   \ . \]
\noindent We consider the Hilbert $C_b(Y)$-bimodule $\wH \otimes_X C_b(Y)$ (defined as in (\ref{def_ea})), and the Hilbert $C_b(Y)$-bimodule $C_b ( Y , p_* \mH  )$ of bounded continuous sections of $p_* \mH$. Then, a natural map $i : \wH \otimes_X C_b(Y) \ra  C_b (Y , p_* \mH)$ is defined, by assigning $[i ( \psi \otimes f)] (y) := ( f (y) \cdot \psi ( p (y) ) , y )$, $y \in Y$. It is clear that $i$ is an injective $C_0(X)$-bimodule map. Furthermore, for every $(v, y) \in p_* \mH$ there is $\psi \in \wH$ such that $v = \psi (\pi (y))$. Thus, we established an inclusion of Hilbert $C_0(Y)$-bimodules
\begin{equation}
\label{pullback1}
\wH \otimes_X C_b(Y)  \hra  C_b (Y , p_* \mH)  \ .
\end{equation}
\noindent In the case in which $X$, $Y$ are compact there is a unital inclusion $C(X) \hra C(Y)$; by assuming $\wH$ full, we find that $i ( \wH \otimes_X C(Y) )$ is full as a Hilbert $C(Y)$-bimodule, and defines a total subset for the continuous field associated with $p_* \mH$. By Rem.\ref{full_bimod} $i$ is surjective, and there is an isomorphism
\begin{equation}
\label{pullback}
\wH \otimes_X C(Y) \simeq  p_* \wH \ ,
\end{equation}
\noindent where $p_* \wH$ denotes the module of continuous sections of $p_* \mH$.

\begin{lem}
Let $\mH \ra X$ be a Hilbert bundle. If $p : Y \ra X$ is a surjective map of locally compact Hausdorff spaces, then there is a monomorphism $C_b(Y) \otimes_X \coh \hra C_b (Y, \wa \mO_{p_* \mH})$. If $X, Y$ are compact (and $\wH$ is full), then $\mO_{p_* \mH}$ is isomorphic to $C(Y) \otimes_X \coh$.
\end{lem}

\begin{proof}
It follows from (\ref{pullback1}) that for every $r,s \in \bN$ there is a monomorphism $C_b(Y) \otimes_X C_0 (X , \mH^{r,s}) \hra C_b ( Y , p_* \mH^{r,s})$. Since there is an obvious isomorphism $C_b(Y) \otimes_X \coh \simeq \mO_{\wH \otimes_X C_b(Y)}$, the lemma follows by applying the universality of the Pimsner algebra to the monomorphism (\ref{pullback1}). The proof in the compact case goes through the same lines, by using (\ref{pullback}).
\end{proof}

\bigskip

We now describe the Toeplitz extension for $\coh$ in terms of \sC algebra bundles. For this purpose, let us make use of the original construction of the Pimsner algebra in \cite{Pim93}; we consider the Fock bimodule
\begin{equation}
\label{fock}
F(\mH) := C_0(X) \oplus \bigoplus \limits_{k = 1}^\infty \wH^k \ .
\end{equation}
\noindent $F(\mH)$ defines in a natural way a continuous field of Hilbert spaces over $X$, with fibres the (separable) Fock spaces $F(\bC^{d(x)}) := \bC \oplus \sum_k^\oplus \bC^{{d(x)}^k}$. 

\noindent Let $LF(\mH)$ denote the \sC algebra of $C_0(X)$-endomorphisms of $F(\mH)$. For every $\psi \in \wH$, we consider the shift operators $\ovl \psi \in LF(\mH)$
\[  \ovl \psi (\psi')  := \psi \otimes \psi' \ \ ,  \]
\noindent $\psi' \in {\wH}^k \subset F(\mH)$. The Toeplitz algebra $\toh$ is defined as the \sC subalgebra of $LF(\mH)$ generated by the operators $\ovl \psi$. Let us denote by $KF(\mH)$ the ideal of compact operators in $LF (\mH)$; by \cite[Thm.3.13]{Pim93}, there is an exact sequence of \sC algebra bundles
\begin{equation}
\label{toe_gen}
0 \ra KF(\mH) \ra \toh \ra \coh \ra 0 \ \ .
\end{equation}
Let us now suppose that $X$ is paracompact. Since $KF(\mH)$ is the \sC algebra of compact operators on the continuous field of Hilbert spaces $F(\mH)$, the Dixmier-Douady class of $KF(\mH)$ vanishes in $H^3(X,\bZ)$ (see \cite[10.7.15]{Dix}). If $\wH$ is countably generated, then $F(\mH)$, $KF(\mH)$ are separable as continuous fields, and $KF(\mH) \otimes \mK$ is locally trivial (\cite[\S 10 p.69]{Dup}). Thus by \cite[10.8.4]{Dix} we find $KF(\mH) \otimes \mK \simeq C_0(X) \otimes \mK$, and obtain the exact sequence
\begin{equation}
\label{toe_gen1}
0 \ra C_0(X) \otimes \mK \ra \toh \otimes \mK \ra \coh \otimes \mK \ra 0 \ \ .
\end{equation}

\noindent Let now in particular $\mE \ra X$ be a vector bundle with an open trivializing cover $\left\{ U_k \right\}$. The local charts $\pi_k : \left. \mE \right|_{U_k} \ra U_k \times \bC^d$ induce by functoriality local charts

\[  \pi_k^F : \left. F(\mE) \right|_{U_k} \ra C_b(U_k) \otimes F(\bC^d)  \ .  \]

\noindent Thus $KF(\mE)$ is locally trivial; by applying \cite[10.8.4]{Dix} we find that $KF(\mE)$ is trivial, and obtain the exact sequence
\begin{equation}
\label{toe_coe}
0 \ra C_0(X) \otimes \mK \ra \toe \ra \coe \ra 0  \ \ .
\end{equation}

Since by \cite[Thm.4.4]{Pim93} $\toh$ is $K$-equivalent to $C_0(X)$, it follows from (\ref{toe_gen1}) that there is an exact sequence $K^0(X) \ra K^0(X) \ra K_0(\coh)$. It is clear that the above $K$-theoretical exact sequence can be regarded as a part of those in \cite[Thm.4.9]{Pim93}.

\bigskip

The DR-algebra $\doh$ is endowed with a canonical shift endomorphism $\widetilde \sigma$, in the following way: if $t \in (\hrs) \subset \doh$, then  
\begin{equation}
\label{def_sigma}
\widetilde \sigma (t) := 1 \otimes t \ \ ,
\end{equation}
\noindent where $1$ is the identity on $\mE$. Note that $\widetilde \sigma$ is well defined over $\doh$ because the left $C_0(X)$-module action coincides with the right one; in fact, in general the tensor product $1 \otimes t$ does not make sense in the setting of Hilbert bimodule operators. Note that by definition $\widetilde \sigma (f) = f$ for every $f \in C_0(X)$; furthermore, it is clear that the canonical endomorphism commutes with the circle action, i.e. $\widetilde \sigma \circ \wa z = \wa z \circ \widetilde \sigma$, $z \in \bT$, so that $\widetilde \sigma (\doh^k) = \doh^k$, $k \in \bZ$. $\widetilde \sigma$ is 'local' in the following sense: if $U$ is any open (closed) subset of $X$ and $\widetilde \sigma_U$ is the canonical endomorphism over ${\widetilde \mO}_{ \left. \mH \right|_U }$, then $\widetilde \sigma_U (t_U) = (\widetilde \sigma (t))_U$ for every $t \in \doh$, where $t_U$ is the restriction of $t$ over $U$ as a vector field.

\bigskip

\begin{rem}
\label{coh_sigma_stable}
The Pimsner algebra $\coh$ is $\widetilde \sigma$-stable. In fact, if $t \in \mK (\hrs) \equiv C_0 (X , \mH^{r,s}) \subset \coh$, then $t = f t'$, $f \in C_0(X)$, $t' \in \mK(\hrs)$ (by \cite[Prop.1.8]{Bla96}), and $\widetilde \sigma (t) = 1 \otimes t = f \otimes t' \in \mK(\hrs)$. We denote by $\sigma$ the restriction of the canonical endomorphism on $\coh$.
\end{rem}

\bigskip

\begin{rem}
\label{end_uni}
It is well-known that the space ${\bf end}_1 \mO_d$ of the unital endomorphisms of $\mO_d$ is homeomorphic to the unitary group $U\mO_d$. Let $X$ be a locally compact, paracompact Hausdorff space, $\mE \ra X$ a vector bundle. If $u \in \doe$ is unitary, then $u \cdot \wE$ is a Hilbert bimodule in $\coe$ with support $1$, thus by universality a $C_0(X)$-endomorphism $\rho_u$ is induced on $\coe$. Viceversa, let us consider a set $\left\{ \psi_i \right\}$ of generators for $\wE$ as in Lemma \ref{generators}. If $\rho$ is a $C_0(X)$-endomorphism of $\coe$, then the net $\left\{ \sum_i \rho (\psi_i) \psi^*_i \right\}$ converges in the strict topology to a unitary of $\doe$. Thus, we established (for $X$ paracompact) a homeomorphism between the space of $C_0(X)$-endomorphisms of $\coe$ and the unitary group of $\doe$, endowed with the usual topologies.
\end{rem}

\bigskip

\begin{rem}
\label{rem_theta_1}
Let us now consider the following operators in $\doh$, which will be used in the sequel:
\begin{equation}
\label{def_theta}
\left\{ 
\begin{array}{l}
\theta_{0,k} := 1 \\

\theta_{1,k} \in (\mH^{k+1} , \mH^{k+1}) \ : \  
\theta (\psi \otimes \psi') := \psi' \otimes \psi \ , \ 
\psi \in \wH , \psi' \in \wH^k \\ 

\theta_{r,k} := \theta_{1,k} \sigma (\theta_{1,k}) \ \cdots \ \sigma^{r-1} (\theta_{1,k})
\in (\mH^{k+r} , \mH^{k+r})

\end{array}
\right.
\end{equation}
\noindent Elementary computations show that every $\theta_{r,k}$ is a symmetry, i.e. $\theta_{r,k} = \theta^*_{r,k} = \theta^{-1}_{r,k}$. Furthermore, if $\psi \in \wH$, then $\theta_{1,k} \cdot (\psi \otimes 1_k) \cdot \psi' = \psi' \otimes \psi =  (1^k \otimes \psi) \cdot \psi'$ for every $\psi' \in \wH^k$ (here $1_k \in (\mH^k , \mH^k)$ denotes the identity map); thus, we get the equality
\begin{equation}
\label{sigma_psi}
\theta_{1,k} \psi = \sigma^k (\psi)  \ \ , \ \ \psi \in \wH
\end{equation}
\noindent Let now $a \in \mK (\mH^r , \mH^r) \subset \coh^0$; then $a$ is a norm limit of sums of terms of the type 
\begin{equation}
\label{decomp}
\psi_{l_1} \cdots \psi_{l_r} \psi_{m_r}^* \cdots \psi_{m_1}^*
\end{equation}
\noindent where $\psi_{l_i} , \psi_{m_j} \in \wH$, $i,j = 1 , \cdots , r$. Elementary computations, performed by applying (\ref{def_theta},\ref{sigma_psi},\ref{decomp}), show that
\begin{equation}
\label{thetark}
\left\{ 
\begin{array}{l}
\mathrm{ad} \theta_{r,k} (a) :=  \theta_{r,k} a \theta_{r,k}^* = \sigma^k (a) \ \ , \ \
             a \in \mK (\mH^r , \mH^r)  \\
\mathrm{ad} \theta_{s,k} (a) = \sigma^k (a) \ , \ s \geq r \ \ .
\end{array}
\right.
\end{equation}
\end{rem}

\bigskip

The canonical endomorphism has been studied in the case of the Cuntz algebra $\mO_d$ (i.e., $X = \left\{ x \right\}$) from the $K$-theoretical viewpoint (see \cite{Cun81}), while certain symmetry properties (related to the operators $\theta_{r,k}$) have been investigated in \cite{DR87,DR89A} for applications to group duality. In \cite[\S 4]{Vas03}, such symmetry properties are studied in the general case $X \neq \left\{ x \right\}$: by using the terminology of \cite{DR89A}, $\widetilde \sigma \in {\bf end}\doh$ has permutation symmetry, but in general it does not have a well defined dimension, and it does not satisfy the special conjugate property unless $\mH$ is a vector bundle with vanishing first Chern class.

\bigskip

Let now $\mE \ra X$ be a vector bundle over a locally compact, paracompact Hausdorff space $X$, and $\mE^\otimes$ denote the category with objects the tensor powers $\mE^r$, $r \in \bN$, and arrows $(\ers)$. It is well known that $\mE^\otimes$ is a tensor category, when equipped with the internal tensor product of vector bundles (see \cite{Ati,Kar}). We claim that the following equality holds in $\doe$:
\begin{equation}
\label{amenable}
(\ers) = (\widetilde \sigma^r , \widetilde \sigma^s) := 
\left\{ t \in \doe : t \widetilde \sigma^r (t') = \widetilde \sigma^s (t') t \ , \ t' \in \doe \right\}  \ ;
\end{equation}
\noindent in other terms, $\mE$ satisfies the {\em amenability property} in the sense of \cite[\S 5]{LR97}. In order to prove (\ref{amenable}), we note that if $t \in (\ers)$, $t' \in (\mE^{r'} , \mE^{s'})$ then obviously $t \otimes t' = 1_s \otimes t' \cdot t \otimes 1_{r'} = t \otimes 1_{s'} \cdot 1_r \otimes t'$, where $1_r$ is the identity over $\mE^r$. Thus, the relation $t \widetilde \sigma^r (t') = \widetilde \sigma^s (t') t$ holds in $\doe$, i.e. $(\ers) \subseteq (\widetilde \sigma^r , \widetilde \sigma^s)$. Viceversa if $t \in (\widetilde \sigma^r , \widetilde \sigma^s)$, we pick sets of generators $\left\{ \psi_{i , h} \right\}$, $\left\{ \varphi_{i , k} \right\}$ respectively for $\wE^s$, $\wE^r$ as by Lemma \ref{generators}, and note that $t$ is the limit (in the strict topology) of the net $\left\{ t_i := \sum_{h,k}  \psi_{i , h} t_{i , hk} \varphi_{i , k}^* \right\}$, where $t_{i , hk} := \psi_{i , h}^* t \varphi_{i , k}$. Since $t_{i , hk} \in (\ii) = C_b(X)$, we find $t_i \in ( \ers )$, so that $t \in (\ers)$. In the same way, by recalling Lemma \ref{K-bundle}, the equality
\begin{equation}
\label{amenable1}
\mK (\ers) = (\sigma^r , \sigma^s) := 
\left\{ t \in \coe : t \sigma^r (t') = \sigma^s (t') t \ , \ t' \in \coe \right\} 
\end{equation}
\noindent is proved.

\bigskip

\begin{prop}
\label{sigma_inner}
Let $X$ be a locally compact, paracompact space, $\mE \ra X$ a vector bundle. Then, the following properties hold:

\begin{itemize}

\item  The canonical endomorphism $\sigma$ of $\coe$ is topologically inner (in the sense of Def.\ref{inner}).

\item  If $(\iota , \mE) \simeq C_b (X , \mE)$ is finitely generated as a Hilbert $C_b(X)$-module, then $\sigma$ is inner (in the sense of Def.\ref{inner}). There is a vector bundle
\[
\beta \mE \ra \beta X \ \ : \ \ \left. \beta \mE \right|_X = \mE \ \ ,
\]
\noindent with natural isomorphisms of Hilbert $C_b(X)$-bimodules $M_\sigma \simeq C_b(X , \mE) \simeq \widehat {\beta \mE}$. Furthermore, there is an isomorphism $\doe \simeq \mO_{\beta \mE}$.

\end{itemize}
\end{prop}

\begin{proof}
It follows from (\ref{amenable}) that $M_\sigma = (\iota ,  \widetilde \sigma) = (\iota , \mE) = C_b (X , \mE)$, so that it is clear that $ \sigma$ is topologically inner. Furthermore, $\sigma$ is inner if and only if $C_b (X , \mE)$ is finitely generated as a Hilbert $C(\beta X)$-module, and $\beta \mE$ is defined by the Serre-Swan theorem. The fact that $\left. \beta \mE \right|_X = \mE$ follows by considering the immersion $\wE \hra C_b (X , \mE) \simeq \widehat {\beta \mE} \simeq \mK (\iota , \beta \mE)$, and the identification of $X$ as a subset of $\beta X$. Now, {\em a priori} $(\beta \mE^r , \beta \mE^s) = \mK (\beta \mE^r , \beta \mE^s) \subseteq C_b (X , \mE^{r,s})$; but it is clear that $\mK (\beta \mE^r , \beta \mE^s)$ is a total subset for $C_b (X , \mE^{r,s})$, so that $\mK (\beta \mE^r , \beta \mE^s) = C_b (X , \mE^{r,s})$. Thus, $(\ers) = C_b (X , \mE^{r,s}) = (\beta \mE^r , \beta \mE^s )$, and $\doe$ is isomorphic to $\mO_{\beta \mE}$.
\end{proof}

\begin{ex}
Let $\mE \ra X$ be a rank $d$ vector bundle. Suppose there is a compact set $Y \subset X$ such that $\left. \mE \right|_{X - Y}$ is trivial; then $C_b (X , \mE)$ is finitely generated.
\end{ex}

Let $\sigma$ be inner in $\coe$, and $\left\{ \psi_l \right\}$ a (finite) set of generators for $p_* \wE$. From the previous proposition, and by definition of inner endomorphism, we obtain the identity
\begin{equation}
\label{sigma_exp}
 \sigma (t) = \sum_l \psi_l t \psi_l^* \ \ , \ \ t \in \doe \ .
\end{equation}

\begin{cor}
Let $p \in (\mE^r , \mE^r)$, $r \in \bN$ be a projection, $\sigma_p (t) := p \otimes t$, $t \in \mK (\ers)$, the shift endomorphism on $\coe$ induced by $p$. Then, $\sigma_p$ is topologically inner; $\sigma_p$ is inner if $p \cdot ( \iota , \mE^r )$ is finitely generated as a Hilbert $C_b(X)$-module.
\end{cor}

\begin{proof}
The result is obtained in the same way as in the previous proposition, by considering the Hilbert $C_0(X)$-bimodule in $\coe$ given by $M_p := p \cdot (\iota , \mE^r)$. It is clear that $\sigma_p$ is induced by $M_p$. Note that $M_p$ has support $p$.
\end{proof}

\begin{rem}
Let $\toe$ denote the Toeplitz algebra of $\mE$. Then, $\wE$ is contained in $\toe$ as a Hilbert $C_0(X)$-bimodule with support $P := \sum_l \psi_l \psi_l^* \neq 1$, so that a topologically inner endomorphism $\tau \in {\bf end} \toe$ is induced. The explicit expression for $\tau$ is analogous to (\ref{sigma_exp}).
\end{rem}

\bigskip

The canonical endomorphism of the Pimsner (Toeplitz) algebra of a vector bundle plays a universal role for inner endomorphisms in the sense of Def.\ref{inner}; the following result is a consequence of Lemma \ref{lem_bima} (with the subsequent remark), and the definition of inner endomorphism:

\begin{prop}
\label{sigma_uni}
Let $\mA$ be a \sC algebra with $X$ denoting the (compact) spectrum of $ZM(\mA)$, $\rho$ an inner endomorphism on $\mA$. Then, $\rho$ extends to an inner endomorphism $\ovl \rho$ over $M(\mA)$. There is a vector bundle $\mE \ra X$, with a morphism of $C(X)$-dynamical systems $(\toe , \tau ) \ra ( M(\mA) , \ovl \rho)$. If $\ovl \rho$ is unital (i.e., $M_\rho$ has support $1$), then there is a $C(X)$-monomorphism $(\coe ,  \sigma ) \hra ( M(\mA) , \ovl \rho)$.
\end{prop}

\bigskip

In order for applications in the next session, we prove some technical results about the Hilbert $\coe^0$-bimodule $\coe^1$ considered in Sec.\ref{pimsner} and Rem.\ref{hilb_bundle}. In the rest of the present section, $X$ will denote a locally compact Hausdorff space, $\mE \ra X$ a vector bundle. Note that $\coe^1$ is generated as a right $\coe^0$-module by every set of generators for $\wE$.

Let us consider the following Hilbert bimodule structure on $\coe^1$:
\[  
a,t \mapsto \sigma^k (a) t  \ \ , \ \  
t,a \mapsto ta \ \ , \ \ 
\left \langle  t , t' \right \rangle := t^* t'
\]
\noindent where $a \in \coe^0$, $t , t' \in \coe^1$, $k \in \bN$, and $\sigma$ denotes as usual the canonical endomorphism on $\coe$. We denote by $\coe^1(\sigma^k)$ the corresponding Hilbert $\coe^0$-bimodule. 

\begin{rem}
\label{coe1_sigma}
In particular, there is an isomorphism of Hilbert $\coe^0$-bimodules $\beta : \wE \otimes_X \coe^0 \ra \coe^1(\sigma)$, $\beta (\psi \otimes a) := \psi a$, $\psi \in \wE$, $a \in \coe^0$, where $\wE \otimes_X \coe^0$ is defined by (\ref{def_ea}). Note in fact that $\beta (\psi \otimes (a'a)) = \psi a'a = \sigma (a') \psi a$ for every $a \in \coe^0$.
\end{rem}

\bigskip

We recall the reader to the notation (\ref{not_tp}). Note that if $\mA$ is a \sC algebra and $\mM$ a Hilbert $\mA$-$I\mA$-bimodule, then for every $s \in [0,1]$ there is a natural covariant epimorphism $\pi_s : \mM \ra \pi_s (\mM)$, induced by the evaluation morphism $ev_s : I\mA \ra \mA$, $ev_s (a) := a(s)$, $a \in I\mA$. Every $\pi_s (\mM)$ is a Hilbert $\mA$-bimodule in the natural way.

\begin{lem}
\label{lem_coe1_homotopy}
Let $\mE \ra X$ be a vector bundle. Then, for every $k \in \bN$ the bimodule $\coe^1$ is homotopic to $\coe^1(\sigma^k)$; i.e., there exists a Hilbert $\coe^0$-$I \coe^0$-bimodule $\mM_k$ such that $\pi_0 (\mM_k) \simeq \coe^1$, $\pi_1 (\mM_k) \simeq \coe^1(\sigma^k)$.
\end{lem}

\begin{proof}
Let $\theta_{r,k} \in \doe$ be the symmetry defined in Rem.\ref{rem_theta_1}. Then, the identity (\ref{sigma_psi}) implies that $\theta_{r,k}$ is the unitary associated with the endomorphism $\sigma^k$ of $\coe$. Now, there is a homotopy $(u_s)_s$, $s \in [0,1]$ in the unitary group of $\doe^0$, connecting $\theta_{r,k}$ with $1$ (see \cite[Prop.4.2.7]{WO}). We denote by 
\[
(\rho_s)_s \ : \ \rho_s (\psi) := u_s \psi \ \ , \ \ \psi \in \wE
\]
\noindent the corresponding homotopy in the space of endomorphisms of $\coe$; since $\rho_s$ is induced by a unitary of $\doe^0$, we find that $\coe^0$ is $\rho_s$-stable for every $s \in [0,1]$. By construction, we have $\rho_0 = id$, $\rho_1 = \sigma^k$. We then consider the space $\mM_k := C ([0,1] , \coe^1)$ of norm continuous maps from $[0,1]$ into $\coe^1$, endowed with the Hilbert $\coe^0$-$I \coe^0$-bimodule structure
\[
\left\{
\begin{array}{ll}

a , \psi \mapsto  \left\{ s \mapsto \rho_s (a) \psi (s) \right\}    \\

\psi , b \mapsto  \left\{ s \mapsto \psi (s) b (s) \right\}    \\

\left \langle \psi , \psi' \right \rangle := 
\left\{ s \mapsto \psi^*(s) \psi'(s) \right\}

\end{array}
\right.
\]
\noindent $\psi , \psi' \in \mM_k$, $a \in \coe^0$, $b \in I \coe^0$, $s \in I$. It is clear that $\mM_k$ satisfied the required conditions.
\end{proof}

\begin{cor}
\label{cor_co1d_hom}
There is a Hilbert $X \mO_d^0$-$IX \mO_d^0$-bimodule connecting $C_0(X , \mO_d^1)$ with the free Hilbert $X \mO_d^0$-bimodule $\oplus^d (X \mO_d^0)$.
\end{cor}

\begin{proof}
We apply the previous lemma to the trivial bundle $\eps (d) := X \times \bC^d$, so that there is a homotopy connetting $\mO_{\eps (d)}^1 = C_0(X , \mO_d^1)$ and $\mO_{\eps (d)}^1 (\sigma)$. Let us prove that $\mO_{\eps (d)}^1 (\sigma)$ can be identified with $\oplus^d (X \mO_d^0)$. At this purpose, recall that $\mO_{\eps (d)}^1 (\sigma)$ coincides with $C_0(X , \mO_d^1)$ as a right Hilbert $(X \mO_d^0)$-module. Now, if $t \in \mO_{\eps (d)}^1 (\sigma)$ and $\left\{ \psi_k \right\}_{k=1}^d \subset \mO_d^1$ is the set of orthogonal isometries generating the Cuntz algebra, then the continuous map $t_k := \left\{ x \mapsto \psi_k^* t (x) \in \mO_d^0 \right\}$ vanishes at infinity, and $t_k \in X \mO_d^0$. Let $a \in X \mO_d^0$; since $\psi_k^* \sigma (a) = a \psi_k^*$, we find $(\sigma (a)t)_k = a t_k$. The desired isomorphism $\beta : \mO_{\eps (d)}^1 (\sigma) \ra \oplus^d (X \mO_d^0)$ is given by $\beta (t) := (t_k)_{k=1}^d$; note in fact that $\beta (ta) = (t_k)_k a$, $\beta (\sigma (a)t ) = a (t_k)_k$, $\left \langle \beta(t) , \beta(t') \right \rangle = \sum_k t^* \psi_k \psi_k^* t' = t^* t'$, where $a \in X \mO_d^0$, $t , t' \in  \mO_{\eps (d)}^1 (\sigma)$.
\end{proof}

\bigskip

We conclude the present section with some considerations in the case in which the base space $X$ is compact. If $\left\{ \psi_l \right\}_l$ is any finite set of generators for $\wE$, then the following relations are satisfied in $\coe$ (see \cite[\S 3]{Pim93}):
\begin{equation}
\label{genrel}
f \psi_l = \psi_l f  \ \ , \ \
\psi_l^*  \psi_m = \left \langle \psi_l , \psi_m \right \rangle \ \ , \ 	\
\sum \limits_l \psi_l  \psi_l^* = 1 \ .
\end{equation}
\noindent Here $f \in C(X)$ and $\left\langle \cdot , \cdot \right\rangle$ denotes the $C(X)$-valued scalar product on $\wE$. By removing the third of (\ref{genrel}), we obtain relations satisfied by the Toeplitz algebra $\toe$. These relations identify $\coe$ (resp. $\toe$) in terms of universal properties, as for Cuntz-Krieger \sC algebras. Note that the matrix $\left \langle \psi_l , \psi_m \right \rangle_{l,m} \in C(X) \otimes \bM_n$ is a projection, actually the one corresponding to the projective module $\wE$ according to the Serre-Swan theorem. Thus we obtain the following result, which is an application of the universality of the Pimsner algebra (\cite[\S 3]{Pim93}, also see \cite[\S 4]{Vas}):

\begin{thm}
\label{thm_genrel}
Let $X$ be a compact Hausdorff space, $E := (E_{l,m}) \in \bM_n \otimes C(X)$ a projection. Then there exists a unique \sC norm over the *-algebra generated by $C(X)$ and a finite set $\left\{  \psi_l \right\}_{l = 1}^n$, with relations
\begin{equation}
\label{genrel1}
\left\{ 
\begin{array}{l}
f \psi_l = \psi_l  f  \\
\psi_l^* \psi_m = E_{l,m} \\
\sum \limits_l \psi_l  \psi_l^* = 1
\end{array}
\right.
\end{equation}
\noindent extending the $\sup$-norm over $C(X)$, and such that the circle action
\[  
\wa z (\psi) := z \cdot \psi_l \ \ , \ \ z \in \bT \ ,
\]
\noindent extends to a $\bT$-action by automorphisms. The \sC algebra so generated is isomorphic to the Pimsner algebra $\coe$, where $\mE \ra X$ is the vector bundle corresponding to $E$ via the Serre-Swan theorem. By removing the last of (\ref{genrel1}), the Toeplitz algebra $\toe$ is obtained.
\end{thm}

The following result is a consequence of the projectivity of $\wE$.

\begin{prop}
Let $X$ be a compact Hausdorff space, $\mE \ra X$ a vector bundle. If $\mE'$ is a vector subbundle of $\mE$, then there is a unital inclusion $\mT_{\mE'} \hra \coe$ of \sC algebra bundles. In particular, there is $n \in \bN$ with a unital inclusion $\toe \hra C(X) \otimes \mO_n$ of \sC algebra bundles.
\end{prop}

\begin{proof}
The first assertion follows from the fact that $\wE'$ is a finitely generated, projective submodule of $\wE$. So that, $\wE' \subset (\iota , \mE)$ is contained in $\coe$ as a finitely generated Hilbert $C(X)$-bimodule with support $P \neq 1$. Thus, the assertion follows from the argument of Lemma \ref{lem_bima} (and the successive remark). The second assertion follows from the fact that $\wE$ is a submodule of the free bimodule $C(X) \otimes \bC^n$, for some $n \in \bN$.
\end{proof}

\begin{rem}
With the notation of Thm.\ref{thm_genrel}, an amplimorphism associated with $\coe^1$ is given by $\phi (a) := (\psi_l^* a \psi_m)_{l,m} \in \bM_n \otimes \coe^0$, $a \in \coe^0$. In the same way, an amplimorphism associated with $\coe^1(\sigma^k)$, $k \in \bN$, is given by $\phi_k (a) := {\mathrm {diag}}_n (\sigma^{k-1} (a)) \cdot E \in \bM_n \otimes \coe^0$, $a \in \coe^0$, where ${\mathrm {diag}}_n (a)$ is the diagonal matrix in $\bM_n \otimes \coe^0$. These facts can be easily checked by considering the map of right Hilbert $\coe^0$-modules $\beta (t) := (\psi_l^* t)_l \in \oplus^n \coe^0$, $t \in \coe^1$, and by evaluating $\beta$ over $a t$, $\sigma^k(a)t$, $a \in \coe^0$.
\end{rem}

\section{Graded $\mO_d$-bundles.}
\label{homotopy}

Aim of the present section is to assign $KK$-theoretical invariants to graded $\mO_d$-bundles, for applications to the case of \sC algebras of vector bundles.

\bigskip

As we remarked in Sec.\ref{intro}, the circle action defined on the Cuntz algebra is full, in the sense that $\mO_d$ is generated as a \sC algebra by the disjoint union of the spectral subspaces $\mO_d^k$, $k \in \bZ$. Furthermore, $\mO_d$ is also semi-saturated in the terminology of \cite{AEE95}, i.e. $\mO_d$ is generated as a \sC algebra by $\mO_d^0$, $\mO_d^1$. Thus $\mO_d$ is endowed with a $\bZ$-grading, and each $\mO_d^k$ can be regarded as a (full) Hilbert bimodule over the UHF-algebra $\mO_d^0$.

We denote as usual by $X$ a locally compact Hausdorff space. If $\mF$ is a $C_0(X)$-algebra, we denote by ${\bf aut}_X \mF$ the group of $C_0(X)$-automorphisms of $\mF$.

\begin{defn}
Let $(\mA , G , \alpha )$ be a \sC dynamical system. An $\mA$-bundle $(\mF , (\pi_x : \mF \ra \mA)_{x \in X} )$ has a {\bf global $G$-action} if there exists an action $\alpha^X : G \ra {\bf aut}_X \mF$, such that $\pi_x \circ \alpha^X = \alpha \circ \pi_x$ for every $x \in X$.
\end{defn}

\begin{ex}
Let $\mE \ra X$ be a $G$-vector bundle, where $G$ is a compact group acting trivially on the base space $X$. Then, every fibre $\mE_x \simeq \bC^d$ is a $G$-module in a natural way. A $G$-action by automorphisms is induced on each fibre $\mO_d$ of $\coe$, and $\coe$ has a global $G$-action. In particular, by considering the action $\bT \times \mE \ra \mE$ induced by the scalar multiplication, we recover the natural circle action (\ref{def_circ_act}) on $\coe$, which is globally defined on $\coe$ as an $\mO_d$-bundle in the sense of the previous definition.
\end{ex}

Let $\mF$ be a  $\mO_d$-bundle carrying a global circle action. Then, $\mF$ is a $\bZ$-graded \sC algebra: for every $k \in \bZ$ we define the spectral subspace $\mF^k := \left\{ t \in \mF : \wa z (t) = z^k t , z \in \bT \right\}$. As for the case considered in Rem.\ref{hilb_bundle}, every $\mF^k$ has a natural structure of $C_0(X)$-Hilbert $\mF^0$-bimodule. Note that the *-algebra $^0 \mF := \cup_k \mF^k$ is dense in $\mF$. We will say in such a case that $\mF$ is a {\em graded} $\mO_d$-bundle. About the next result, see also \cite[Thm.3.1]{AEE95}.

\begin{prop}
\label{thm_f1}
Let $X$ be a locally compact, paracompact Hausdorff space, $\mF$ a graded $\mO_d$-bundle over $X$. Then $\mF$ is graded $C_0(X)$-isomorphic to the Pimsner algebra $\mO_{\mF^1}$. Furthermore, if $\mB$ is a graded $\mO_d$-bundle, then $\mF$ is graded isomorphic to $\mB$ if and only if the Hilbert bimodules $\mF^1$, $\mB^1$ are covariantly isomorphic.
\end{prop}

\begin{proof}
About the first point, note that $\mF$ admits a locally finite, compact cover $\left\{ U_i \right\}$, with local charts $\pi_i : \mF_i := \left. \mF \right|_{U_i} \ra C(U_i) \otimes \mO_d$. Thus, $\mF_i^1$ is finitely generated as a right $\mF_i^0$-module by a set of orthogonal isometries $\left\{ \psi_k^i \right\}_{k=1}^d \subset \mF^1_i$. By using a partition of unity, we obtain a set $\left\{ \psi_l \right\}$ of generators for $\mF^1$, with $\psi_l^* \psi_m \in \mF^0$: thus, the net $\sum_l \psi_l \psi_l^*$ converges to the identity of the multiplier algebra of $\mF$. In the same way, if $r \in \bN$, by defining $\psi_L := \psi_{l_1} \cdots \psi_{l_r} \in \mF^r$ we obtain that $\mF^r$ is generated by the set $\left\{ \psi_L \right\}$: in fact, $\sum_L \psi_L \psi_L^* \ra 1$, and if $t \in \mF^r$ then $\sum_L \psi_L (\psi_L^* t) \ra t$, with $(\psi_L^* t) \in \mF^0$.

Thus, $\mF^1$ is contained in $\mF$ as a full Hilbert $\mF^0$-bimodule with support $1$; by universality, there is a graded $C_0(X)$-monomorphism $\lambda : \mO_{\mF^1} \ra \mF$. Furthermore, $\lambda$ is surjective, since if $t \in \mF^r$, $r \in \bN$, then by the above considerations $t$ belongs to $\lambda (\mO_{\mF^1})$; the same is true for $\mF^{-r} = (\mF^r)^*$.

About the second point, note that if $\tau : \mF \ra \mB$ is a graded isomorphism, then it is obvious that $(\tau_1 , \tau_0) : \mF^1 \ra \mB^1$ is a covariant isomorphism, where $\tau_k := \left. \tau \right|_{\mF^k}$, $k = 0,1$. Viceversa, if $\mF^1$, $\mB^1$ are covariantly isomorphic, then it follows from the previous point and Lemma \ref{lem_funct} that $\mF$ is graded isomorphic to $\mB$.
\end{proof}

\begin{cor}
\label{morita}
With the above notations, $\mF^r$ is an imprimitivity $\mF^0$-$\mF^0$-bimodule for every $r \in \bN$.
\end{cor}

\begin{proof}
It suffices to consider the identifications $\mK (\mF^r , \mF^r) \simeq \mF^r \cdot (\mF^r)^* \simeq \mF^0$.
\end{proof}

Note that we also established that the circle action over $\mF$ is semi-saturated. As a consequence of the previous result, we have a version of \cite[Thm.2.5]{Pim93}:

\begin{cor}
Let $\mE \ra X$ be a vector bundle. Then there is a $C_0(X)$-isomorphism $(\coe , \bT )  \simeq ( \mO_{\coe^1} , \bT )$.
\end{cor}

From Prop.\ref{thm_f1}, we have an interpretation of the set of isomorphism classes of graded $\mO_d$-bundles in terms of covariant isomorphism classes of Hilbert bimodules. Thus, a description in terms of $KK$-groups becomes natural. 

As a first step, we consider the zero grade algebra. If $\mA$ is a \sC algebra, we denote by $H^1(X , {\bf aut}\mA)$ the set of isomorphism classes of $\mA$-bundles over $X$ (see for example \cite[Thm.2.1]{Vas} for a justification of such a notation). $H^1(X , {\bf aut}\mA)$ has a distinguished element, called $0$, corresponding to the trivial $\mA$-bundle.

If $\mF$ is a graded $\mO_d$-bundle over $X$, we define
\begin{equation}
\label{def_delta_0}
\delta_0 (\mF) := [ \mF^0 \otimes \mK ] \in H^1(X,{\bf aut} ( \mO_d^0 \otimes \mK) ) \ \ .
\end{equation}
\noindent Thus, the equality $\delta_0 (\mF) = \delta_0 (\mB)$ is intended in the sense that $\mF^0$, $\mB^0$ are stably isomorphic as $\mO_d^0$-bundles.

\bigskip

\begin{rem}
Let $X$ be a pointed, compact, connected $CW$-complex such that the pair $(X , x_0)$, $x_0 \in X$, is a homotopy-cogroup. We denote by $SX$ the (reduced) suspension. In order for more compact notations, we define $X^\bullet := X - \left\{ x_0 \right\}$. It follows from a result by Nistor (\cite[\S 5]{Nis92}) that 
\[
H^1 (SX , {\bf aut}\mO_d^0) \simeq   [X , {\bf aut} \mO_d^0 ]_{x_0}  
                            \simeq   KK_1 ( C \mO_d^0 \ , \ X^\bullet \mO_d^0 ) \ ,
\]
\noindent where $C \mO_d^0 := \left\{ ( z , a) \in \bC \oplus C_0 ( [0,1) , \mO_d^0 ) \ : \ a (0) = z 1 \right\}$ is the mapping cone. We also find
\[
[X , {\bf aut} (\mO_d^0 \otimes \mK) ]_{x_0}  
\simeq   KK_0 ( \mO_d^0 \ , \ X^\bullet \mO_d^0 ) \ .
\]
\noindent In particular, when $X$ is the $(i-1)$-sphere, we get $H^1 ( S^i , {\bf aut} \mO_d^0  ) = \left\{ 0 \right\}$, as proved also in \cite[Thm.1.15]{Tho87}.
\end{rem}

\bigskip

In the sequel, we will use the concepts and the notations introduced in Sec.\ref{kk_groups}.

\bigskip

Let $X$ be a $\sigma$-compact metrisable space; then, every graded $\mO_d$-bundle $\mF$ is separable and $\sigma$-unital, and $\mF^1$ is countably generated as a Hilbert $\mF^0$-bimodule (in fact, $\mF$ is finitely generated over compact subsets). Furthermore by Cor.\ref{morita} $\mF^0$ acts on the left over $\mF^1$ by elements of $\mK ( \mF^1 , \mF^1  )$. Thus, we define
\begin{equation}
\label{def_delta_1}
\delta_1 (\mF) := [( \mF^1 ,0 ) ] \in KK ( X ; \mF^0 , \mF^0 ) \ \ .
\end{equation}
\noindent With an abuse of notation, in the sequel we will denote by $\delta_1 (\mF)$ also the class of $(\mF^1 , 0)$ in $KK_0 ( \mF^0 , \mF^0 )$ obtained by {\em forgetting} the $C_0(X)$-structure. Note that since $\mF^1$ is an imprimitivity bimodule, we find that $\delta_1 (\mF)$ is invertible; thus, the Kasparov product by $\delta_1 (\mF)$ defines an automorphism on $KK_0 ( \mA , \mF^0 )$ for every \sC algebra $\mA$. In particular, $\delta_1 (\mF) \in {\bf aut} K_0 (\mF^0)$.

\bigskip

From \cite[Thm.4.9]{Pim93} and Prop.\ref{thm_f1}, we get an exact sequence for the $KK$-theory of $\mF$. It is clear that in the case in which $\mF$ is the Pimsner algebra of a vector bundle, we can directly apply \cite[Thm.4.9]{Pim93} by replacing $\mF^0$ with $C_0(X)$.

\begin{cor}
\label{cor_ex_seq}
For every separable \sC algebra $\mA$, and graded $\mO_d$-bundle $\mF$, the following exact sequence holds:
\begin{equation}
\label{ex_seq}
\xymatrix{
           KK_0 (\mA , \mF^0)
		    \ar[r]^-{{1 - \delta_1(\mF)}}
		 &  KK_0 (\mA , \mF^0)
		    \ar[r]^-{i_0}
		 &  KK_0 (\mA , \mF)
		    \ar[d]^-{\delta_0}
		 \\ KK_1 (\mA , \mF)
		    \ar[u]^-{\delta_1}
		 &  KK_1 (\mA , \mF^0)
		    \ar[l]_-{i_1 }
		 &  KK_1 (\mA , \mF^0)                             
		    \ar[l]_{{1- \delta_1(\mF)}}
}
\end{equation}
\noindent where $i_*$ are the morphisms induced by the inclusion $\mF^0 \hra \mF$, and $\delta_*$ are the connecting maps induced by the $KK$-equivalence between $\mF^0$, $\mT_{\mF^1}$.
\end{cor}

\bigskip

We introduce a notation. Let $\mF$, $\mB$ be graded $\mO_d$-bundles with $\delta_0(\mF) = \delta_0(\mB)$; then, there is a $C_0(X)$-algebra isomorphism $\alpha : \mF^0 \otimes \mK \ra \mB^0 \otimes \mK$, and a ring isomorphism $\alpha_* : KK (X ; \mF^0 , \mF^0) \ra KK (X ; \mB^0 , \mB^0)$ is defined. We write

\begin{equation}
\label{def_iso_delta}
\delta (\mF) = \delta (\mB)  \ \
\Leftrightarrow \ \ 
\delta_0 (\mF) = \delta_0 (\mB) \ \
{\mathrm {and}} \ \
\alpha_* \delta_1(\mF) = \delta_1 (\mB) \ ;
\end{equation}

\noindent note that we used the stability of $KK(X ; - , - )$, so that we identified $[\mF^1] \in KK (X ; \mF^0 , \mF^0)$ with $[\mF^1 \otimes \mK] \in KK (X ; \mF^0 \otimes \mK , \mF^0 \otimes \mK )$. The tensor product of $\mF^1$ by $\mK$ is intended as the external tensor product of Hilbert bimodules. Note that $\alpha_* \delta_1(\mF) = [ (\mF^1 \otimes \mK)_\alpha ]$, where $(\mF^1 \otimes \mK)_\alpha$ is the pullback bimodule defined as in Ex.\ref{ex_pullback}.

\bigskip

\noindent Let us now denote by 
\begin{equation}
\label{def_smi}
i : RK^0(X) \ra KK (X ; \coe^0 , \coe^0)
\end{equation}
\noindent the structure morphism (\ref{def_i}).

\begin{thm}
\label{delta1_coe}
Let $\mE \ra X$ be a rank $d$ vector bundle. Then 
\begin{equation}
\label{eq_main}
\delta_1 (\coe) = i[\mE]
\end{equation}
\noindent Furthermore, if $\mF$ is the trivial $\mO_d$-bundle, then $\delta_0 (\mF) = 0$, $\delta_1 (\mF) = [d]$.
\end{thm}

\begin{proof}
It follows from Lemma \ref{lem_coe1_homotopy} that $\coe^1$ is homotopic to $\coe^1(\sigma)$. The left $\coe^0$-module action over the Hilbert $\coe^0$-$I \coe^0$-bimodule $\mM := C ([0,1] , \coe^1)$ realizing the homotopy is by elements of $\mK ( \mM , \mM )$, thus $\mM$ defines a Kasparov cycle in $KK (X ; \coe^0 , I \coe^0)$. So that, $[\coe^1] = [\coe^1(\sigma)] \in KK(X ; \coe^0 , \coe^0)$. On the other side, it follows from Rem.\ref{coe1_sigma} that $\coe^1 (\sigma)$ is isomorphic to the Hilbert $\coe^0$-bimodule $\wE \otimes_X \coe^0$, defining the class $i [\mE] \in KK (X; \coe^0 , \coe^0)$. Finally, if $\mF \simeq X \mO_d$, then the same argument as above shows that $\mF^1 \simeq C_0 (X , \mO_d^1)$ is homotopic to $\oplus^d (X \mO_d^0)$ (see Cor.\ref{cor_co1d_hom}). This last bimodule has class $[d]$ in $KK ( X ; X \mO_d^0 , X \mO_d^0 )$, with the notation (\ref{eq_nkk}).
\end{proof}

We emphasize the fact that the r.h.s. of (\ref{eq_main}) is a $KK$-theoretical class defined independently by the notion of Pimsner algebra (see \S \ref{kk_groups}).

\bigskip

\begin{rem}
\label{rem_kcoe0}
Let $X$ be paracompact. The $K$-theory of $\coe^0$ can be obtained in the following way: as a first step, note that for every $r \in \bN$ we find $\mK (\mE^r , \mE^r) \otimes \mK \simeq C_0(X) \otimes \mK$ (in fact, $\mK (\mE^r , \mE^r) \otimes \mK$ is a $\mK$-bundle with trivial Dixmier-Douady class). Thus we have the inductive limit
\[  
K_0 (\coe^0) = K^0(X) \stackrel{i_*}{\ra} 
               K^0(X) \stackrel{i_*}{\ra} 
               K^0(X) \stackrel{i_*}{\ra} \cdots
\]
\noindent where $i_*$ is the adjoint map of the inclusion $i : \mK (\mE^r , \mE^r) \hra \mK (\mE^{r+1} , \mE^{r+1})$ (Sec.\ref{pimsner}). We denote the elements of $K_0(\coe^0)$ by 
\[
(y_r) \equiv (y_r)_{r \in \bN} \in K_0(\coe^0) \ \ ,
\]
\noindent $y_r \in K^0(X) \simeq K_0(\mK (\mE^r , \mE^r))$, $r \in \bN$ (so that, $y_{r+1} = i_* y_r$ for every $r \geq$ of a given $r_0 \in \bN$). Let us now suppose for simplicity that $X$ is compact; then, the automorphism on $K_0 (\coe^0)$ induced by $i[\mE] = \delta_1 (\coe)$ is given by
\[  
i[\mE] \ (y_r) := ([\mE] \otimes y_r ) \ \ .
\]
\noindent In particular, if $1_r \in K_0(\mE^r , \mE^r)$ is the class of the identity of $(\mE^r , \mE^r)$, then the image of $1_r$ w.r.t. the isomorphism $K_0(\mE^r , \mE^r) \simeq K^0(X)$ is the class $1 \in K^0(X)$ corresponding to the trivial line bundle. We denote by $(1)$ the corresponding sequence in $K_0(\coe^0)$. Thus, we get a map
\[
K_0 (X) \ni [\mE]   \mapsto   i[\mE] \ (1) = ([\mE]) \in K_0(\coe^0) \ \ .
\]
\end{rem}

\bigskip
\bigskip

Let $\mF$ be a graded $\mO_d$-bundle. It is clear that in general the elements of $KK ( X ; \mF^0 , \mF^0 )$ do not arise from grade-one components of graded $\mO_d$-bundles. Anyway, they can be recognized by considering any open trivializing cover $\left\{ U \right\}$ for $\mF$, and by noting that the conditions $\delta_0 ( \mF^0 |_U ) = 0$, $\delta_1 (\mF^1 |_U ) = [d]$ hold, thanks to the previous theorem. On the other side, at a global level note that every graded $\mO_d$-bundle defines an imprimitivity bimodule over the zero-grade algebra (Cor.\ref{morita}). 

\bigskip

Let us now consider the stabilized bundle $\mF \otimes \mK$; then the $\bT$-action extends in the natural way $t \otimes a \mapsto \wa z (t \otimes a) := z (t) \otimes a$, $t \in \mF$, $z \in \bT$, $a \in \mK$. By Fourier analysis, it is clear that
\[
(\mF \otimes \mK)^k := 
\left\{ \ovl t \in \mF \otimes \mK : \wa z (\ovl t) := z^k \ovl t \right\} =
\mF^k \otimes \mK \ ,
\]
\noindent $k \in \bZ$. In particular, $\mF^1 \otimes \mK$ is a $C_0(X)$-Hilbert $(\mF^0 \otimes \mK)$-bimodule in $\mF \otimes \mK$. Now, the \sC subalgebra generated by $\mF^1 \otimes \mK$ coincides with $\mF \otimes \mK$; furthermore, it follows from the argument of Ex.\ref{ex_son} that $\mF^1 \otimes \mK$ has support $1$. By universality of the Pimsner algebra, we conclude that there is a $C_0(X)$-isomorphism $( \mF \otimes \mK , \bT  ) \simeq ( \mO_{\mF^1 \otimes \mK} , \bT )$. It is also clear that $\mF^1 \otimes \mK$ is an imprimitivity $C_0(X)$-Hilbert $(\mF^0 \otimes \mK)$-bimodule. By Lemma \ref{pic_pim}, we obtain:

\begin{prop}
\label{prop_out}
Let $\mF$, $\mB$ be graded $\mO_d$-bundles over a $\sigma$-compact Hausdorff space $X$. Suppose there is a $C_0(X)$-isomorphism $\alpha : \mF^0 \otimes \mK \ra \mB^0 \otimes \mK$ such that the pullback bimodule $(\mF^1 \otimes \mK)_\alpha$ is outer conjugate to $\mB^1 \otimes \mK$ in ${\bf {Pic}}_X (\mB^0)$. Then there is a $C_0(X)$-isomorphism $(\mF \otimes \mK , \bT ) \ra (\mB \otimes \mK , \bT )$.
\end{prop}

\begin{cor}
\label{cor_fr}
Suppose there is $r \in \bN$, $r \geq 1$, with isomorphisms of Hilbert $(\mB^0 \otimes \mK )$-bimodules $\mB^r \otimes \mK \simeq (\mF^r \otimes \mK)_\alpha $, $\mB^{r+1} \otimes \mK \simeq (\mF^{r+1} \otimes \mK)_\alpha$. Then there is a $C_0(X)$-isomorphism $(\mF \otimes \mK , \bT ) \ra (\mB \otimes \mK , \bT )$.
\end{cor}

\begin{proof}
In order for economize in notations, we write $\widetilde \mF^r := \mF^r \otimes \mK$. It follows from Cor.\ref{morita} that $\widetilde \mB^r $, $\widetilde \mB^{r+1} $ (resp. ${\widetilde \mF}^r_\alpha$, ${\widetilde \mF}^{r+1}_\alpha$) are imprimitivity $\widetilde \mB^0$-bimodules. Now, for every $r,s \in \bN$ there is a natural isomorphism $\widetilde \mF^{r+s}  = \widetilde \mF^r \otimes_{\widetilde \mF^0 } \widetilde \mF^s$ of Hilbert $\widetilde \mF^0$-bimodules, and the analogue holds for spectral subspaces in $\mB$. Thus,
\[
\widetilde \mB^1 \otimes_{\widetilde \mB^0 } \widetilde \mB^r  \simeq 
\widetilde \mB^{r+1}  \simeq 
{\widetilde \mF}^{r+1}_\alpha \simeq
{\widetilde \mF}^{r}_\alpha \otimes_{\widetilde \mB^0 } {\widetilde \mF}^1_\alpha \ ,
\]
\noindent and it suffices to apply the previous proposition (note that $\left\{ \widetilde \mB^r  \right\} = \left\{ {\widetilde \mF}^r_\alpha \right\}$).
\end{proof}

\begin{rem}
\label{rem_f}
Let $\theta(\mF) := \theta_{\mF^1 \otimes \mK}$ be the automorphism of $\mF^0 \otimes \mK$ defined according to (\ref{iso_pic}). By recalling the construction (\ref{def_malfa}), and by universality of the Pimsner algebra, we obtain the isomorphism

\[ 
( \mF \otimes \mK , \bT )
\simeq 
( ( \mF^0 \otimes \mK )\rtimes_{\theta(\mF)} \bZ  , \bT ) \ . 
\]

\noindent Thus the Picard group of the stabilized zero-grade bundle, modulo outer conjugacy classes, describes the graded isomorphism classes of stabilized $\mO_d$-bundles (Prop.\ref{prop_out}): the map
\[
{\bf {Pic}}_X (\mF^0)  \ra  KK ( X ; \mF^0 , \mF^0 )  \ \ , \ \ 
\theta (\mF) \mapsto \delta_1 (\mF)
\]
\noindent gives a measure of the accuracy of the class $\delta_1$ in describing the graded isomorphism classes of stabilized $\mO_d$-bundles.
\end{rem}

\bigskip

We now prove some results in the case of Pimsner algebras of vector bundles.

\begin{prop}
\label{prop_k}
Let $\mE \ra X$ be a rank $d$ vector bundle. If $\mE' \ra X$ is a rank $d$ vector bundle with $\delta (\coe) = \delta (\mO_{\mE'})$, then $i[\mE] = i[\mE']$. In particular, $i[\mE] = i[\mE']$ if there is a $C_0(X)$-isomorphism $( \coe \otimes \mK , \bT ) \ra ( \mO_{\mE'} \otimes \mK , \bT )$. If the structure morphism (\ref{def_smi}) is injective, then $[\mE] = [\mE'] \in RK^0 (X)$.
\end{prop}

\begin{proof}
We denote by $i' : RK^0(X) \ra KK( X ; \mO_{\mE'}^0 , \mO_{\mE'}^0 )$ the structure morphism (\ref{def_i}). Let $\alpha_*$ be the ring isomorphism 

\[  
\alpha_* :  KK (X; \coe^0 , \coe^0) \ra KK(X; \mO_{\mE'}^0 , \mO_{\mE'}^0)
\]

\noindent induced by an isomorphism between $\coe^0$, $\mO_{\mE'}^0$. Then, by (\ref{knature}) we find $i'[\mE'] = \delta_1 (\mO_{\mE'}) = \alpha_* \delta_1 (\coe) = \alpha_* i[\mE] = i'[\mE]$. Since $\alpha_*$ is injective, we obtain $i[\mE] = i[\mE']$.
\end{proof}

\begin{rem}
\label{rem_keri}
From the previous proposition it is clear that the computation of the group $\ker i \subseteq RK^0(X)$ is of interest. A way to check the injectivity of $i$ is to regard at $i[\mE]$ as an automorphism of $K_0(\coe^0)$: in fact, if $X$ is compact we find $i [\mE] \ (1) := ([\mE]) \in K_0(\coe^0)$ (see Rem.\ref{rem_kcoe0}). If the map $K^0 (X) \ni [\mE] \mapsto ([\mE]) \in K_0(\coe^0)$ is injective, we obtain
\[
i[\mE] = i [\mE'] \ \Rightarrow \
([\mE]) = ([\mE']) \ \Rightarrow \
[\mE] = [\mE'] \in K^0(X) \ .
\]
\end{rem}

\begin{rem}
Note that if $X = \bullet$ reduces to a single point, then $\mE \simeq \bC^d$, $RK^0(\bullet) = \bZ$, and $KK ( \bullet , \coe^0 , \coe^0 ) = KK_0 ( \mO_d^0 , \mO_d^0  ) = {\bf end} \  \bZ \left[ \frac{1}{d}  \right]$. Thus we find the structure morphism $i : \bZ \ra {\bf end} \  \bZ \left[ \frac{1}{d}  \right]$, where $i (k)$ is the multiplication by $k$ in $ \bZ \left[ \frac{1}{d}  \right]$ (see \cite[10.11.8]{Bla},\cite{Cun81}). In particular, $i[\mE]$ is the multiplication by $d$ (that is an automorphism on $\bZ \left[ \frac{1}{d}  \right]$). More generally, let $X$ be compact and suppose that $\coe^0 \simeq X \mO_d^0$ (for example, this is always true if $X$ is an $i$-sphere, see \cite[Thm.1.15]{Tho87}). If $K^0(X)$ is finitely generated, we can apply the Kunneth Theorem \cite[\S 23]{Bla} and obtain $KK_0 ( \coe^0 , \coe^0 ) = {\bf end} \left( K^0(X) \otimes \bZ \left[ \frac{1}{d}  \right] \right)$: so that $i[\mE]$ is the multiplication by the class $[\mE] \otimes 1 \in K^0(X) \otimes \bZ \left[ \frac{1}{d}  \right]$.
\end{rem}

\begin{prop}
\label{prop_cw}
Let $X$ be a (compact) finite dimensional $CW$-complex, $\mE , \mE' \ra X$ rank $d$ vector bundles with $[\mE] = [\mE'] \in K^0(X)$. Then, $\delta (\coe) = \delta (\mO_{\mE'})$, and there is a $C(X)$-isomorphism $(\coe \otimes \mK , \bT ) \ra ( \mO_{\mE'} \otimes \mK , \bT )$.
\end{prop}

\begin{proof}
For every $r \in \bN$, we have $[\mE^r] = [\mE'^r] \in K^0(X)$. $X$ being a finite dimensional $CW$-complex, we find that $\mE^r$ is isomorphic to $\mE'^r$ for every $r$ greater than a fixed $r_0 \in \bN$ (\cite[VIII.Thm.1.5]{Hus}). Thus, by Rem.\ref{zero_grade}, we obtain that there is an isomorphism $\alpha : \coe^0 \ra \mO_{\mE'}^0$. Now, $\delta_1 (\mO_{\mE'}) = i' [\mE'] = i' [\mE] =  \alpha_* i[\mE] = \alpha_* \delta_1 (\coe)$, having applied (\ref{knature}). Finally, the fact that $\coe$, $\mO_{\mE'}$ are stably isomorphic follows by Cor.\ref{cor_fr}.
\end{proof}

\section{Final remarks (Work in Progress).}
\label{final}

It is of interest to perform a classification of the \sC algebras $\coe$, up to stable isomorphism. For this purpose two important points, that will be approached in a successive work (\cite{Vasf}), are the following:

\begin{enumerate}

\item the computation of $\ker i \subseteq RK^0(X)$, that gives a characterization in terms of representable $K$-theory of vector bundles having stably graded isomorphic Pimsner algebras (recall Rem.\ref{rem_kcoe0}, Rem.\ref{rem_keri}).

\item the computation of the representable $KK$-group $KK (X ; \coe , \coe)$, in analogy with the recent work by Kirchberg and Blanchard (\cite{Kir00,BK02}). It is clear that the main tools remain the exact sequences \cite[Thm.4.9]{Pim93}, which have to be generalized to the $KK(X ; - , -)$-bifunctor.

\end{enumerate} 

A related problem is whether the condition $\delta (\mF) = \delta (\mB)$ implies a graded stable isomorphism for the graded $\mO_d$-bundles $\mF$, $\mB$ (see Rem.\ref{rem_f}).

\bigskip

It has been remarked in Sec.\ref{basics} that the canonical endomorphism over $\doe$ plays a sort of universal role w.r.t. the notion of inner endomorphism introduced in Def.\ref{inner}. 
This fact is well known, in the particular case in which a given endomorphism $\rho$ is implemented by a set of orthogonal isometries, i.e. -- with the terminology of Def.\ref{inner} -- $M_\rho$ is a free Hilbert bimodule (see \cite{DR89A} and related references). The corresponding crossed products involve the Cuntz algebras, and play an important role in the Doplicher-Roberts duality for compact groups (see \cite{DR87,DR89,DR89A}). 
The \sC algebras $\coe$ play an analogous role, in the more general case in which there is a nontrivial centre and $M_\rho$ is not free. The corresponding duality characterizes tensor categories arising from non-compact group actions over Hilbert bimodules. About these aspects, see \cite{Vas03} and related references.

\bigskip

An other field of application for the Pimsner algebra of a vector bundle is in the context of projective multiresolution analysis in the sense of \cite{PR03}. In fact, with methods analogous to \cite{BJ98}, it is possible to construct representations of Pimsner algebras associated with certain vector bundles over the torus $\bT^n := \left. \bR^n \right/ \bZ^n$, arising from a projective multiresolution analysis. This topic will be developed elsewhere.

\bigskip
\bigskip

{\bf Acknowledgments.} The core of the present work started as a degree thesis supervised by S. Doplicher; the author would like to thank him for the initial idea, support and encouragements. The author also thanks E. Blanchard, A. Valette, A. Kishimoto for fruitful email correspondences, and P. Goldstein for discussions and comments during the period he was PostDoc in Rome.


\end{document}